\documentclass[10pt,reqno]{amsart}
\usepackage[margin=1.1in]{geometry}
\usepackage{amsmath,amssymb,amsthm,array}
\usepackage[T1]{fontenc}
\usepackage{hyperref}

\numberwithin{equation}{section}
\newtheorem{theorem}{Theorem}[section]
\newtheorem{proposition}[theorem]{Proposition}
\newtheorem{lemma}[theorem]{Lemma}
\newtheorem{corollary}[theorem]{Corollary}
\theoremstyle{definition}
\newtheorem{definition}[theorem]{Definition}

\newtheorem{remark}[theorem]{Remark}
\newtheorem{example}[theorem]{Example}

\newcommand{\Z}{\mathbb{Z}}
\newcommand{\C}{\mathbb{C}}
\newcommand{\Rep}{\operatorname{Rep}}
\newcommand{\Aut}{\operatorname{Aut}}
\newcommand{\id}{\mathrm{id}}
\newcommand{\ev}{{\bar 0}}
\newcommand{\od}{{\bar 1}}
\newcommand{\Vev}{V_{\ev}}
\newcommand{\Vod}{V_{\od}}
\newcommand{\M}{\mathcal{M}}
\newcommand{\R}{\mathcal{R}}
\newcommand{\Q}{\mathcal{Q}}
\newcommand{\fus}{\boxtimes}
\newcommand{\ol}[1]{\overline{#1}}

\title[Twisted Verlinde formulas for cyclic orbifolds and the order-two simple current]
{Twisted Verlinde formulas for cyclic orbifolds and the order-two simple current of 
a vertex operator superalgebra}
\author{A. Zuevsky}
\date{}
\address{Institute of Mathematics, Czech Academy of Sciences, \v{Z}itn\'{a} 25, Prague, Czech
		Republic}
\email{zuevsky@yahoo.com}

\begin{document}

\begin{abstract}

Let $V=\Vev\oplus J$ be a simple vertex operator superalgebra with $\Vev$ simple, rational 
and $C_2$-cofinite, and let $G$ be a finite cyclic group acting faithfully on $\Vev$.  
For arbitrary twisted $V$-modules we prove: 
(1) for $G$ of arbitrary order, a twisted Verlinde formula for $\Vev$; 
(2) unitarity of the orbifold $S$-matrix block in generating sectors;
(3) as a consequence, a formula for $N^{X_3}_{X_1,X_2}\pm N^{J\boxtimes X_3}_{X_1,X_2}$ 
as sums over $J$-orbits of $G$-stable untwisted modules, in which $J$-fixed points enter 
with weight two and in the minus sign, with a parity selection rule 
(read in $\mathbb Z/2\mathbb Z$) and a dimension sum rule;
(4) simplicity of $V^G$.
These results concern the $G$-graded fusion category of $\Vev$ together with the order two  
invertible object $J$ with $\theta_J=-\dim_qJ$. 
  The statements are tested numerically on six 
 examples (odd and even orders up to $15$, both signs of $\dim_qJ$, $J$-fixed modules). 
\end{abstract}



\maketitle


\section{Introduction}

The Verlinde formula expresses the fusion rules of a rational vertex operator algebra through 
the modular $S$-matrix \cite{V, H1}. This paper concerns two extensions. 

\emph{Twisted modules.} For a finite automorphism group $G$ of a rational, $C_2$-cofinite 
vertex operator algebra, Dong and Lin \cite{DL} proved, for $G$ of prime order,  
a twisted Verlinde formula for fusion rules among $g$-twisted modules, in terms of the 
orbifold $S$-matrix and a sum over $G$-stable untwisted modules.

\emph{Superalgebras.} For $V=\Vev\oplus\Vod$ with regular even part, Creutzig, Kanade and 
McRae \cite{CKM} realized $V$ as a superalgebra object (an extension by an order two simple 
current $J=\Vod$) in the modular tensor category of $\Vev$-modules and obtained Verlinde 
formulas for dimensions $N^+$ and superdimensions $N^-$ of spaces of intertwining operators
 among untwisted and parity-twisted $V$-modules.

In this paper, we combine the two: 
 we extend the even-part formula to cyclic groups of arbitrary order and 
derive the super case formula, including $J$-fixed points, for the combined situation.

\subsection{Previous works}

The relation between the $S$-matrix of $V^G$ and the $S$-matrix of the orbifold theory was
 determined in \cite{DRX1} for finite groups and used in 
\cite{DL} for prime order. Categorically, Verlinde formulas for $G$-crossed modular fusion 
categories are due to Turaev \cite{T} (see also \cite{BBCW}) and Deshpande-Mukhopadhyay 
\cite{DeM}.  
  It was shown in \cite{Mc, DL} that the category of twisted $\Vev$-modules 
is $G$-crossed modular. 
 Theorem \ref{thmTV} below is the vertex algebraic realization  
of the $G$-crossed Verlinde formula for cyclic $G$. 

Our Lemma \ref{lemorbrel} and Theorem \ref{thmTV} are derived  from modular invariance of 
twisted traces \cite{DLM3} and the classification of $V^G$-modules 
\cite{DY, MT, DRX}. For $|G|$ prime, Lemma \ref{lemorbrel} (b) 
 (with identical conjugation conventions, see Remark \ref{remconv}), Theorem \ref{thmTV},
 and the case distinction in part 3 of its proof corresponds to \cite{DL}. 
    
 For composite order new proven topics are:
 the treatment of arbitrary stabilizers $G_X$ in the orbifold relation, 
the restriction to stable columns by character orthogonality in the stabilizers, 
and the scope of the individual coefficient formula (see Remark \ref{remTVopen}). 
 Lemma \ref{lemorbrel} is connected to the relations of \cite{DRX1} for general finite $G$.  
We give a self-contained derivation which keeps 
track of the stabilizers. For permutation orbifolds see \cite{DXY}.

\subsection{Main results of this paper}

Let $G$ be cyclic of order $n$, acting faithfully on $\Vev$ 
(assumptions of Subsection \ref{secsetting}). 
Let $\Q$ be the set of $G$-stable irreducible untwisted $\Vev$-modules, $S_{X,Q}$ the orbifold 
$S$-matrix ($X$ twisted, $Q\in\Q$), and $d=\dim_qJ=\pm1$. We then have the following results. 

\smallskip
\textit{Even part, arbitrary cyclic order:} 

\textbf{Theorem \ref{thmTV}.} For $g_1$, $g_2$  
generating $G$, $g_3=g_1g_2$ and $M_i\in\M(g_i)$, the average of $N^{M_3}_{N_1,N_2}$ over 
$N_i$ in the $G$-orbits of $M_i$ equals 
$\sum_{Q\in\Q}S_{M_1,Q}S_{M_2,Q}\ol{S_{M_3,Q}}/S_{V_\ev,Q}$. If two of the three modules are 
$G$-stable (automatic for $n$ a prime power, or if all modules are stable) this is the 
individual fusion coefficient. 

\textbf{Proposition \ref{propunit}:} 
for $g$ generating $G$, the block $(S_{X,Q})_{X\in\M(g),Q\in\Q}$ is square and unitary.

\smallskip
\textit{A consequence of Theorem \ref{thmTV} and the  
$J$-equivariance of Section \ref{secsign}:} 

\textbf{Theorem \ref{thmmain} (a).} 
 For $X_i$ as above, with 
$N^\pm:=N^{X_3}_{X_1,X_2}\pm N^{J\boxtimes X_3}_{X_1,X_2}$ and $\pi(X)\in\Z/2\Z$ 
the parity type (see Definition \ref{defsigns}),
\[
N^{\pm}=\delta_{\pi(X_1)+\pi(X_2),\,\pi(X_3)}\Bigl(4\!\!\sum_{[Q]\ \text{free},\ 
\eta_Q=\pm1}\!\!f_0(Q)\;+\;2\!\!\sum_{Q\ \text{fixed},\ \eta_Q=\pm1}\!\!f_0(Q)\Bigr),
\qquad f_0(Q)=\frac{S_{X_1,Q}S_{X_2,Q}\ol{S_{X_3,Q}}}{S_{V_\ev,Q}},
\]
where the equality of parity types is in $\Z/2\Z$ (equivalently 
$\epsilon_{X_1}\epsilon_{X_2}\epsilon_{X_3}=1$). Here $[Q]$ runs over $J$-orbits of 
size two in $\Q$, where ``fixed'' means $J\fus Q\simeq Q$. Fixed points exist only for 
$d=1$ and have $\eta_Q=-1$ (see Lemma \ref{lemfixed}). For $g_X$ generating $G$, $X$ is 
$J$-fixed iff $S_{X,Q}=0$ whenever $\eta_Q=-1$ (Lemma \ref{lemfixed} (iii)).

 \textit{Dimension sum rule:} 

\textbf{Proposition \ref{propsum}} 
and parity selection rule:  \textbf{Corollary \ref{corselect}},
 and \textbf{Lemma \ref{lemfixedsimple}:} $V^G$ is simple whenever $V$ is.

 \textbf{Theorem \ref{thmmain} (b), Corollary \ref{corproduct}:} Under conditions 
\ref{hypdict}, $N^+$ and $N^-$ are the dimension and superdimension of the space of 
intertwining operators among the corresponding $g$-twisted and $gP$-twisted $V$-modules. 
 The fusion-product statement is based in addition on conditions \ref{hypmono}.
 

Theorems \ref{thmTV} and \ref{thmmain} (a), Proposition \ref{propunit} 
and the sum rule use only the assumptions on $\Vev$ and $U=\Vev^G$ and 
on modular invariance of twisted traces.  
   The statements that concern twisted $V$-modules 
 (Theorem \ref{thmmain} (b), Corollary \ref{corproduct}, Proposition \ref{propinner} (iii)) 
  use conditions \ref{hypdict}-\ref{hypL}. 
 For inner automorphisms they reduce to the untwisted theory \cite{CKM, HKL} provided 
Li $\Delta$-operator results 
\cite{Li, DLM1} hold for superalgebras, including the fusion of two twisted modules
(conditions \ref{hypL}).

\subsection{The examples:} 
The numerical evidence is given in Section \ref{secexamples}.  The algebra of 
the proofs and the stated conventions are tested in all examples. The stable-column formula 
of Lemma \ref{lemorbrel} (b) is tested independently of $S^U$ in one family 
(Example \ref{exA}, from twisted traces). The stabilizer dependent parts are tested  
 only for $n=2$ (Example \ref{exB}). The composite order example with a  
proper stabilizer (Example \ref{exF}) is a product and tests only the statements. 
Section \ref{secopen} lists the remaining open problems.

\section{Preliminaries}

\subsection{Vertex operator superalgebras}
\label{secsetting} 

Let $V=\Vev\oplus J$ be a vertex operator superalgebra with parity automorphism $P$. 
 Automorphisms preserve the grading. Ideals are $\mathbb Z/2\mathbb Z$-graded, and ``simple'' 
means that there are no non-trivial graded ideals. Twisted modules are as in \cite{DLM2, X} 
with the convention of \cite{DLM3} for $g\leftrightarrow g^{-1}$.  
$g$-twisted $V$-modules are called \emph{parity-untwisted}, $gP$-twisted ones are called 
\emph{parity-twisted}. In this paper we adopt the following two assumptions 
(in what follows we call them the first/second assumptions of Subsection \ref{secsetting}):  
 
\emph{First:} $V$ is a simple, $\tfrac12\Z$-graded vertex operator superalgebra of CFT type. 
$\Vev$ is simple, rational, $C_2$-cofinite, self-dual, of CFT type.

\emph{Second:} $G=\langle\sigma\rangle\subset\Aut(V)$ is cyclic of order $n\ge2$ and acts 
faithfully on $\Vev$.

Set $U:=\Vev^G$, $J:=\Vod$, $J_0:=\Vod^G$. 
We have the following consequences of two assumptions above:  

\quad 1. $\Vev$ is $g$-rational for every finite order automorphism $g$ 
(with the assumption that $\Vev$ is rational and $C_2$-cofinite,  
which is the first assumption above) \cite{DL}.  

 \quad 2.  For a solvable (here cyclic) subgroup  
$G\subset\Aut(\Vev)$ the fixed-point algebra $U=\Vev^G$ is again simple, rational, 
$C_2$-cofinite, self-dual and of CFT type (with assumption that $\Vev$ is simple, rational, 
$C_2$-cofinite, self-dual, of CFT type, i.e., the first assumption) 
applied to the vertex operator algebra
 $\Vev$), \cite{DL, CM, M}. 

\quad 3.  For $W\in\{\Vev,U\}$ the category $\mathcal C_W$ 
is modular, the $S$-matrix is symmetric and unitary, and the Verlinde formula holds 
(with assumptions as above), \cite{DL, H1, H2, DLN}. 

These results are statements about vertex operator algebras. The superalgebra $V$  
enters through the $\Vev$-module $J$ and the superalgebra object structure of \cite{CKM}
 (see Proposition \ref{factJ}). The following lemma describes simplicity of $V^G$.

\begin{lemma}
\label{lemfixedsimple}

Let $V$ be a simple vertex operator superalgebra and $G$ a finite group of automorphisms. 
Then $V^G=U\oplus J_0$ is a simple vertex operator superalgebra.

\end{lemma}

\begin{proof}

\emph{ 1) Generation:} Let $M$ be a module for a vertex superalgebra $V$, $c\in M$ 
homogeneous. We show that $S:=\mathrm{span}\{u_nc:\ u\in V,\ n\in\Z\}$ is a submodule. 
Let $a$, $b\in V$ be homogeneous, $\varepsilon=(-1)^{|a||b|}$. The Borcherds identity gives, 
for all $m$, $n$, $k\in\Z$, 

\begin{equation}
\label{eqborch}
\sum_{j\ge0}\binom mj\,(a_{n+j}b)_{m+k-j}\,c
=\sum_{j\ge0}(-1)^j\binom nj\Bigl[a_{m+n-j}\,b_{k+j}\,c-(-1)^{n}
\varepsilon\,b_{n+k-j}\,a_{m+j}\,c\Bigr], 
\end{equation}
where signs appear through $\varepsilon$. By the truncation axioms there are $N$ with 
$a_ic=0$ for all $i\ge N$, and $N'$ with $b_ic=0$ for all $i\ge N'$. Also, $a_{n+j}b=0$ for 
$j\gg0$. Take $m\ge N$. Then the second sum on the right of \eqref{eqborch} vanishes, and 
the left side belongs to $S$. Hence, for all $m\ge N$ and all $n$, $k$: 
$\sum_{j=0}^{J}(-1)^j\binom nj\,a_{m+n-j}b_{k+j}c\in S$, where $J=N'-k-1$ 
(the terms with $j>J$ vanish). Fix $r\in\Z$ and $k$, and set $n=r-m$: 
\[
\sum_{j=0}^{J}p_j(m)\,a_{r-j}b_{k+j}c\in S\quad\text{for all }m\ge N,
\qquad p_j(m):=(-1)^j\binom{r-m}{j}.
\]
Each $p_j$ is a polynomial in $m$ of  degree $j$, thus the $(J+1)\times(J+1)$ matrix 
$(p_j(N+i))_{0\le i,j\le J}$ is invertible. Hence, every vector $a_{r-j}b_{k+j}c$ belongs to 
$S$. In particular, for $j=0$, $a_rb_kc\in S$ for all $r$, $k$. 
 Thus, $S$ is stable under all modes of all homogeneous elements. 
It contains $c=\mathbf 1_{-1}c$. This proves 1).

\emph{2) We show that $I'$ is a nonzero graded ideal of $V$:} let $0\ne I\subset V^G$ be a graded 
ideal and $I':=\mathrm{span}\{u_ni:\ u\in V,\ n\in\Z,\ i\in I\}$. Applying 1) to $M=V$  
 as adjoint module, $I'=\sum_{i\in I}\mathrm{span}\{u_ni\}$ is a sum of left ideals. 
 Hence, it is stable under action all $v_m$, $v\in V$. It is spanned by homogeneous 
elements, hence it is graded. It contains $I$ ($u=\mathbf 1$), hence is nonzero. It is 
stable under action of $D=L(-1)=\omega_0$ which is a mode of $\omega\in V$. By skew-symmetry, 
$i'_nv=\pm\sum_{k\ge0}\frac{(-1)^{n+1+k}}{k!}D^k\bigl(v_{n+k}i'\bigr)$ for 
$i'\in I'$, $v\in V$, which belongs to $I'$ since $I'$ is stable under 
$v_m$ and $D$. Thus, $I'$ is a two-sided graded ideal. By simplicity of $V$, $I'=V$.

\emph{3) Averaging and both inclusions:} let $\mathrm{av}=|G|^{-1}\sum_{g\in G}g$. 
For $u\in V$, $i\in I\subset V^G$, $g\in G$: $g(u_ni)=(gu)_n(gi)=(gu)_ni$. Hence,  
$\mathrm{av}(u_ni)=\mathrm{av}(u)_n\,i\in I$, since $\mathrm{av}(u)\in V^G$ and $I$ 
is an ideal of $V^G$. Thus, $\mathrm{av}(I')\subset I$. Since $I'=V$, 
$V^G=\mathrm{av}(V)=\mathrm{av}(I')\subset I$. 
Together with $I\subset V^G$ this gives $I=V^G$. 
\end{proof}

\begin{lemma}
\label{lemGP}

Under the second assumption of Subsection \ref{secsetting}, 
$G\cap\langle P\rangle=\{1\}$. Hence,  
$G\times\langle P\rangle\simeq\Z/n\times\Z/2$ acts faithfully on $V$ and the sectors $g$, 
$gP$ (with $g\in G$) are all distinct.

\end{lemma}

\begin{proof}
$P$ acts trivially on $\Vev$ and $G\to\Aut(\Vev)$ is injective, thus $P\in G$ would lead 
to $P=1$.
\end{proof}
Note that $J$ is a $G$-stable $\Vev$-module ($\sigma|_J:J\to J\circ\sigma$ is an isomorphism).

Let us adopt the following conventions.
 In the braided category $\mathcal C_{\Vev}$ of \cite{H2} the ribbon 
twist of a simple module $A$ is $\theta_A=e^{2\pi ih_A}$ (where $h_A$ the conformal weight),   
$\dim_q$ is the ribbon-categorical dimension, and $c$ is the braiding. For an invertible $J$ 
with $J\otimes J\simeq\mathbf 1$, $c_{J,J}\in\C^\times$ is a scalar, and the ribbon relation
 (curl) gives $\theta_J=c_{J,J}\dim_qJ$ \cite{BK}, since  
$\theta_{J\otimes J}=\theta_{\mathbf 1}=1=\theta_J^2c_{J,J}^2$ and $\dim_qJ=\pm1$, 
$c_{J,J}=\pm1$.
 Then we recall the following proposition. 

\begin{proposition}
\label{factJ}

\begin{enumerate}

\item[(i)]  $J$ is a simple current of order two for $\Vev$. Since  
$V^G=U\oplus J_0$ is a simple vertex operator superalgebra (Lemma \ref{lemfixedsimple})  
with regular even part $U$, $J_0$ is a simple current of order two for $U$ \cite{CKM}. 

\item[(ii)] $\theta_J=-\dim_qJ\in\{\pm1\}$. 

\item[(iii)]  for a simple object $A$ of $\mathcal C_{\Vev}$ or $\mathcal C_U$, 
$\dim_qA=S_{A,\mathbf 0}/S_{\mathbf 0,\mathbf 0}$ \cite{H2, DJX}. 

\end{enumerate}
\end{proposition} 

For $(ii)$, indeed, $V=\Vev\oplus J$ is a 
\emph{superalgebra} object of $\mathcal C_{\Vev}$, \cite{CKM}. The multiplication 
$\mu:J\otimes J\to\Vev$ (nonzero, since $J\otimes J\simeq\Vev$ is simple and $V$ is simple) 
satisfies $\mu\circ c_{J,J}=-\mu$ (supercommutativity on the odd part corresponding to the 
sign $(-1)^{|u||v|}$ in skew-symmetry). Hence, $c_{J,J}=-1$ and $\theta_J=-\dim_qJ$. 
In particular, $\mathcal M_{J,J}=\id$.
 Note that $(iii)$ is the identification between categorical dimensions 
and modular-transformation data.  
 Below $d:=S_{J,V_\ev}/S_{V_\ev,V_\ev}$ is defined through $S$, and (iii) identifies 
it with $\dim_qJ$. It is also used in (ii) to relate $d$ to $\theta_J$. 

\subsection{The categories involved}
\label{seccat}
 
$\mathcal C_{\Vev},\ \mathcal C_U$ are modular tensor categories of ordinary modules \cite{H2}. 
$S^U$ is the (symmetric, unitary) $S$-matrix of $\mathcal C_U$, defined through modular 
transformation of the trace functions $Z_Y(v,\tau)$, $v\in U$.
  $\Rep(\Vev)$ are modules of the commutative algebra object $\Vev\in\mathcal C_U$ \cite{DL}. 
Its simple objects are the irreducible $g$-twisted $\Vev$-modules, $g\in G$.
  $F_0$ is induction $F_0=\Vev\fus_U(\cdot): \mathcal C_U\to\Rep(\Vev)$, left adjoint to the 
forgetful functor.
 $\Rep(V)$, $F_V$ are $\mathbb Z/2\mathbb Z$-graded modules of the superalgebra 
$A=\Vev\oplus J$ in $\Rep(\Vev)$ and induction $F_V(M)=M\oplus J\fus M$. 
Note that it is used under conditions \ref{hypdict}.  

Next, we recall from \cite{DL} (see also \cite{KO, DLXY}) 

\begin{proposition} 
\label{factrep}

Assume that $\Vev$ simple, rational, $C_2$-cofinite, self-dual, 
of CFT type, and $G$ finite abelian group of automorphisms. Then:  

\begin{enumerate}

\item[(a)]  $\Vev$ is a 
commutative algebra in $\mathcal C_U$. $\Rep(\Vev)$ is a tensor category (a fusion category, 
in particular semisimple with finitely many simple objects). Every irreducible $g$-twisted 
$\Vev$-module ($g\in G$) is an object of $\Rep(\Vev)$, two such are equivalent objects iff 
they have the same $g$ and are isomorphic as twisted modules, and every simple object of 
$\Rep(\Vev)$ is an irreducible $g$-twisted $\Vev$-module for some $g$. 

\item[(b)]  the multiplicity $N^{M_3}_{M_1,M_2}$ of $M_3$ in 
$M_1\fus M_2$ equals $\dim I(M_3;M_1,M_2)$, the dimension of the space of intertwining operators of type $\binom{M_3}{M_1\,M_2}$ among twisted modules.

\item[(c)] $N^{M_3}_{M_1,M_2}=0$ unless $g_3=g_1g_2$. Hence,  
$\Rep(\Vev)=\bigoplus_g\Rep(\Vev)_g$ is $G$-graded.

\item[(d)] $N^{M_3\circ h}_{M_1\circ h,M_2\circ h}=N^{M_3}_{M_1,M_2}$ for 
$h\in G$. Consequently, comparing multiplicities in the semisimple category $\Rep(\Vev)$, 
one obtains $(M_1\fus M_2)\circ h\simeq(M_1\circ h)\fus(M_2\circ h)$ as objects.  

\item[(e)] the induction functor $F_0:\mathcal C_U\to\Rep(\Vev)$, $F_0(Y)=\Vev\fus_UY$, is exact and a tensor functor, and $\operatorname{Hom}_{\Rep(\Vev)}(F_0(Y),N)=\operatorname{Hom}_{\mathcal C_U}(Y,N)$ for $N\in\Rep(\Vev)$ regarded as a $U$-module.

\end{enumerate}
\end{proposition}
 
Note that $(d)$ is stated and proved in \cite{DL} for arbitrary $g_1$, $g_2$, $g_3$, $h\in G$. 
The proof is the application of the Jacobi identity and does not use that $|G|$ is prime.

\begin{definition} 
\label{deffusion}

\emph{Fusion coefficients} $N^{M_3}_{M_1,M_2}$ are defined as $\dim I(M_3;M_1,M_2)$,   
 equivalently, by (b), multiplicities in $\Rep(\Vev)$. The tensor product $\fus$ of 
$\Rep(\Vev)$ is used through (a), (e) and for the invertible object $J$. A tensor 
product $\fus_V$ of \emph{twisted $V$-modules} assumes \ref{hypdict}.   

\end{definition}

\subsection{Twisted modules, stabilizers and $U$-modules}
\label{secorb}

Let $\M(g)$ be the irreducible $g$-twisted $\Vev$-modules and $\R=\bigsqcup_{g\in G}\M(g)$.
 For $X\in\M(g)$ let $G_X=\{h\in G:X\circ h\simeq X\}$ (cyclic, and $g \in G_X$ by\cite{DL})  
and $O_X=\{X\circ a\}$ its $G$-orbit, 
$|O_X|=n/|G_X|$. $X$ is \emph{stable} if $G_X=G$. Set $\Q:=\{Q\in\M(1):G_Q=G\}$.

Choose a representation $\varphi_X$ of $G_X$ on $X$ with 
$\varphi_X(h)Y_X(u,z)\varphi_X(h)^{-1}=Y_X(hu,z)$. 
 It is a projective representations of cyclic groups lift,  
$\varphi_X$ is unique up to $\widehat{G_X}$.
 Let $\psi_s(\sigma)=e^{-2\pi is/n}$, thus $\widehat G=\{\psi_s\}$,
 and set $X=\bigoplus_{\chi\in\widehat{G_X}}X_\chi$. 
 For stable $X$ set $X_s=X_{\psi_s}$. For $\Vev$ and $J$ we normalize $\varphi$
 as the action of $\sigma$. Then $(\Vev)_0=U$ and $J_0=\Vod^G$. 
 Recall the following proposition \cite{DRX, DY, MT} (and \cite{DL} for $n$ prime). 

\begin{proposition} 
\label{factclass}

The statements below hold for every $g\in G$ and every $X\in\M(g)$
 with the stabilizer $G_X$: 

\begin{enumerate}

\item[(i)] each $X_\chi$, $\chi\in\widehat{G_X}$, is a nonzero irreducible $U$-module, and $X_\chi\not\simeq X_{\chi'}$ for $\chi\ne\chi'$.

\item[(ii)] for $X\circ a\in O_X$ one has $(X\circ a)_\chi=X_\chi$,  
 with the convention $\varphi_{X\circ a}=\varphi_X$.

\item[(iii)] $X_\chi\simeq Y_\psi$ iff $Y\in O_X$ and $\psi=\chi$. 
For $Y\in O_X$ the groups $G_Y=G_X$ coincide and $\chi$, $\psi$ 
are compared as characters of this group. 

\item[(iv)] every irreducible $U$-module is of the form $X_\chi$.

\end{enumerate}
\end{proposition}

 In \cite{DL} these statements are quoted for $|G|$ prime 
(Theorems 5.2-5.3, where $G_X\in\{1,G\}$). 
 In \cite{DY, MT, DRX} treat finite groups with a projective representation of 
the stabilizer on $X$ in general. 
 For a \emph{cyclic} stabilizer $G_X$ every projective representation lifts to a linear one.
 $H^2$ of a cyclic group with values in $\C^\times$ vanishes.  
 Thus $\varphi_X$ may be taken to be a representation, unique up to $\widehat{G_X}$. 
 We use (i)-(iv) for cyclic $G$ of arbitrary order in this form (see also \cite{DRX, MT}).
  
Concerning the convention in (ii): $X\circ a$ has the same underlying space 
and $Y_{X\circ a}(v,z)=Y_X(av,z)$. For $v\in U$, $Y_{X\circ a}(v,z)=Y_X(v,z)$,  
 thus the $U$-actions agree. For $k\in G_X$, 
 $\varphi_X(k)Y_{X\circ a}(v,z)\varphi_X(k)^{-1}=Y_X(kav,z)=Y_X(akv,z)=Y_{X\circ a}(kv,z)$ 
 for abelian $G$. Thus, $\varphi_X(k)$ satisfies the defining relation of a representation on 
$X\circ a$ with the same stabilizer $G_{X\circ a}=G_X$. Hence, the eigenspaces $X_\chi$ are
 the same subspaces with the same $U$-action, and the labels $\chi$ are preserved. 
By (i) and (iv), a $U$-module $Y$ determines $(O_X,\chi)$ uniquely by (iii). 
 In each orbit there is  one irreducible $U$-module for each $\chi\in\widehat{G_X}$.

\begin{lemma}
\label{lemF0}
$F_0(X_\chi)=\bigoplus_{N\in O_X}N$, each with multiplicity one.
\end{lemma}
\begin{proof}
For a simple $N\in\Rep(\Vev)$ (an irreducible twisted module, in any sector $h\in G$)
 Proposition \ref{factrep} (e) gives 
$\operatorname{Hom}(F_0(X_\chi),N)=\operatorname{Hom}_U(X_\chi,\bigoplus_\psi N_\psi)$. 
 By Proposition \ref{factclass} (i)-(iii) it has dimension one if $N\in O_X$ and $0$ otherwise.
 Then we use semisimplicity.  
\end{proof}

\subsection{The orbifold $S$-matrix}
\label{secS}

We then use the following notations \cite{DL}.
 For commuting $g$, $h\in G$, let $\M(g,h)=\{M\in\M(g):M\circ h\simeq M\}$, 
a set of isomorphism classes of individual modules. 
 For $M\in\M(g,h)$ and $v\in\Vev$ set 
 $Z_M(v,(g,h),\tau)=\operatorname{tr}_Mo(v)\varphi_M(h)q^{L(0)-c/24}$. 
By \cite{DLM3}, the functions $Z_M(v,(g,h),\tau)$ ($v\in\Vev$), 
for all commuting pairs $(g,h)$ and all $M\in\M(g,h)$, are linearly independent, 
and the transformation $\tau\mapsto-1/\tau$ acts by
\begin{equation}\label{eqlaw}
Z_X(v,(g,k),-\tfrac1\tau)=\tau^{\mathrm{wt}[v]}\sum_{L\in\M(k,g^{-1})}S_{[X;(g,k)],[L;(k,g^{-1})]}\,Z_L(v,(k,g^{-1}),\tau),\qquad X\in\M(g,k),
\end{equation}
where $\M(k,g^{-1})=\{L\in\M(k):g\in G_L\}$ and the bracket $[X;(g,k)]$ 
 denotes the module $X$ together with the sector $(g,k)$ and the chosen $\varphi_X(k)$. 
 The coefficients are unique since the functions on the right hand side are linearly independent 
 as functions of $(v,\tau)$ with $v\in\Vev$. 
 For $Q\in\Q$ set $S_{X,Q}:=S_{[X;(g_X,1)],[Q;(1,g_X^{-1})]}$ (row $X$, column $Q$).
 For $g_X=1$ it is the ordinary $S$-matrix of $\Vev$. 
 For a stable $W\in\M(h)$ and $h\in G_X$ set $S^{[h]}_{X,W}:=S_{[X;(g_X,h)],[W;(h,g_X^{-1})]}$.

Let us now describe the restriction to $U$. For $v\in U\subset\Vev$ and 
$a\in G$ one has $Z_{L\circ a}(v,\cdot)=Z_L(v,\cdot)$ (see Proposition \ref{factclass} (ii)). 
Thus, restricted to $v\in U$, the functions of the modules in one $G$-orbit coincide, 
and the sum of the coefficients over an orbit can be read off.
  Lemma \ref{lemorbrel} is a statement about $v\in U$ and uses these orbit sums. 
 
We will use the following conventions.  
  (i) For twisted modules we follow \cite{DLM3}. 
  (ii) $\psi_s(\sigma)=e^{-2\pi is/n}$,
  $X_\chi=\{w:\varphi_X(h)w=\chi(h)w\ \forall h\in G_X\}$. 
  (iii) The transformation $S$ sends the sector $(g,k)$ to $(k,g^{-1})$. 
  All factors $\psi(g^{-1})$ below arise from 
  $Z_L(v,(k,g^{-1}),\tau)=\sum_\psi\psi(g^{-1})Z_{L_\psi}(v,\tau)$, 
  and all factors $\ol{\chi(k)}$ from the projection 
  $Z_{X_\chi}=|G_X|^{-1}\sum_k\ol{\chi(k)}Z_X(\cdot,(g,k),\cdot)$. 
 (iv) $S^U$ is defined by 
 $Z_Y(v,-1/\tau)=\tau^{\mathrm{wt}[v]}\sum_{Y'}S^U_{Y,Y'}Z_{Y'}(v,\tau)$,
 $v\in U$. (v) $T$ acts by $e^{2\pi i(h-c/24)}$. 
 (vi) Three kinds of $S$ occur: 
 the orbifold entries $S_{X,Q}$ (row twisted or untwisted, column $Q\in\Q$), 
 the matrix $S^U$ of $\mathcal C_U$, and the untwisted $S$-matrix of $\Vev$ 
 (the case $g_X=1$ of the first).
 These conventions are confirmed independently in Example \ref{exA}  
 from twisted traces (see Remark \ref{remconv}). Then the following lemma concerning 
 the orbifold relation holds. 

\begin{lemma}
\label{lemorbrel}
Let $X\in\M(g)$, $\chi\in\widehat{G_X}$.
\begin{enumerate}
\item[(a)] If $N\in\M(h)$ and $\psi\in\widehat{G_N}$, then $S^U_{X_\chi,N_\psi}=0$ unless $h\in G_X$ and $g\in G_N$.
\item[(b)] If $W\in\M(h)$ is stable, $h\in G_X$ and $t\in\Z/n$, then
\[S^U_{X_\chi,W_t}=|G_X|^{-1}\,\ol{\chi(h)}\,\psi_t(g^{-1})\,S^{[h]}_{X,W}.\]
In particular, for $Q\in\Q$: \quad $S^U_{X_\chi,Q_t}=|G_X|^{-1}\psi_t(g^{-1})\,S_{X,Q}$.
\item[(c)] If $N\in\M(h)$ has stabilizer $G_N\ni g$, $h\in G_X$, $\psi\in\widehat{G_N}$, then, with the convention $\varphi_L=\varphi_N$ for $L\in O_N$,
\[S^U_{X_\chi,N_\psi}=|G_X|^{-1}\,\ol{\chi(h)}\,\psi(g^{-1})\sum_{L\in O_N}S_{[X;(g,h)],[L;(h,g^{-1})]};\]
(b) is the case $O_N=\{W\}$.
\end{enumerate}
\end{lemma}

\begin{proof}
In what follows let $v\in U$. 

\emph{1) Projection onto $X_\chi$:} for $v\in U$ the operator $o(v)$ commutes with 
$\varphi_X(k)$, $k\in G_X$. The projection onto $X_\chi$ is 
$P_\chi=|G_X|^{-1}\sum_{k\in G_X}\ol{\chi(k)}\varphi_X(k)$. Note that $G_X$ is cyclic, 
thus this is the orthogonal projection for the characters of $G_X$. Hence,
 $Z_{X_\chi}(v,\tau)=|G_X|^{-1}\sum_{k\in G_X}\ol{\chi(k)}Z_X(v,(g,k),\tau)$. 
 This is the origin of the factor $|G_X|^{-1}$ and of $\ol{\chi(k)}$. 
 The sum is over $k\in G_X$, thus $X\in\M(g,k)$ for every term.

\emph{2) Transformation:} apply \eqref{eqlaw} to each term. 
The sector $(g,k)$ becomes $(k,g^{-1})$, and $L$ ranges over $\M(k,g^{-1})$, 
i.e., over $k$-twisted modules with $g\in G_L$.

\emph{3) Decomposition of the transformed traces:} for such $L$, $\varphi_L(g^{-1})$ acts 
on $L_\psi$ by $\psi(g^{-1})$, thus 
$Z_L(v,(k,g^{-1}),\tau)=\sum_{\psi\in\widehat{G_L}}\psi(g^{-1})Z_{L_\psi}(v,\tau)$. 
This is the origin of $\psi(g^{-1})$. Altogether, 
\[
Z_{X_\chi}(v,-\tfrac1\tau)=\tau^{\mathrm{wt}[v]}|G_X|^{-1}\sum_{k\in G_X}\ol{\chi(k)}\sum_{L\in\M(k,g^{-1})}S_{[X;(g,k)],[L;(k,g^{-1})]}\sum_{\psi\in\widehat{G_L}}\psi(g^{-1})Z_{L_\psi}(v,\tau).
\]

\emph{4) Comparison of coefficients:} by definition, 
$Z_{X_\chi}(v,-\frac1\tau)=\tau^{\mathrm{wt}[v]}\sum_YS^U_{X_\chi,Y}Z_Y(v,\tau)$ 
over the irreducible $U$-modules $Y$, and the functions $Z_Y(v,\tau)$, $v\in U$,
 for pairwise inequivalent $Y$ are linearly independent \cite{Z, H1}. 
 By Proposition \ref{factclass} (iii), (iv), the $U$-module $L_\psi$ is equivalent to 
$N_{\psi'}$ iff $L\in O_N$ (hence, $k=h$, since the orbit belongs to one sector) 
and $\psi=\psi'$. Distinct pairs (sector $k$, orbit, and character) 
give inequivalent $U$-modules. 
 Thus, no cancellation between different terms occurs. Consequently,
 the coefficient of $Z_{N_\psi}$ is obtained from the terms with $k=h$ and $L\in O_N$.

\emph{5)} if $h\notin G_X$, no term has $k=h$. Thus, $S^U_{X_\chi,N_\psi}=0$. 
 If $g\notin G_N$ then no $L\in O_N$ belongs to $\M(h,g^{-1})$ (as $G_L=G_N$), 
thus again $0$. This gives (a).
  Otherwise, the terms $L\in O_N$ all equal $Z_{N_\psi}$ on $U$, their coefficients add,
 and we obtain (c). For $O_N=\{W\}$ this gives (b).
\end{proof}

Note that whole orbifold $\rho(S)$ depends on the choices of the representations 
$\varphi_X$ in every sector $(g,k)$ which is proven using 
the existence of the matrix \eqref{eqlaw} in the chosen bases. The following lemma 
describes orbit invariance and dependence on $G$. 

\begin{lemma} 
\label{lemorbinv}
(a) For $X\in\M(g)$, $a\in G$ and $Q\in\Q$: $S_{X\circ a,Q}=S_{X,Q}$.

(b) The entry $S_{X,Q}$ depends only on $X$, on $g=g_X$ and on the operator 
$\varphi_Q(g^{-1})$. In particular, if $H\le G$ contains $g_X$ and $Q$ is $G$-stable, 
the entry calculated with respect to $H$ (with $\varphi_Q|_H$) equals the entry computed 
 with respect to $G$. The set $\Q=\M(1,G)$, however, depends on the group. 
\end{lemma}

\begin{proof}
(a) By Proposition \ref{factclass} (ii), $(X\circ a)_\chi=X_\chi$ and 
$G_{X\circ a}=G_X$. Lemma \ref{lemorbrel} (b) gives 
$|G_X|^{-1}\psi_t(g^{-1})S_{X\circ a,Q}=S^U_{X_\chi,Q_t}=|G_X|^{-1}\psi_t(g^{-1})S_{X,Q}$. 
(b) $S_{X,Q}$ is the coefficient of $Z_Q(v,(1,g^{-1}),\tau)$ 
 in the expansion \eqref{eqlaw} of $Z_X(v,(g,1),-1/\tau)$ over the linearly 
independent functions $Z_L(v,(1,g^{-1}),\tau)$, $L\in\M(1,g^{-1})$. 
 This expansion involves only $g$.
\end{proof}

Let us make a remark concerning two different normalization issues.
  (1) \emph{Orbit convention} for $\varphi_{X\circ a}=\varphi_X$, Proposition \ref{factclass} 
(ii): it affects only the identification of $U$-components inside an orbit and leaves all 
entries $S_{X,Q}$ unchanged (see Lemma \ref{lemorbinv} (a)). 
 (2) \emph{Column normalization:} the entries $S_{X,Q}$ with $g_X\ne1$ depend on 
$\varphi_Q$. If $\varphi_Q$ is replaced by $\lambda(\sigma)\varphi_Q$,
 $\lambda=\psi_a\in\widehat G$, then $Q_t=Q'_{t+a}$. 
 The eigenvalue $\psi_t(\sigma)$ of $\varphi_Q(\sigma)$ becomes $\psi_{t+a}(\sigma)$ 
for $\varphi'_Q(\sigma)$. Comparing Lemma \ref{lemorbrel} (b) for the two labelings gives
\begin{equation}
\label{eqnormchange}
S'_{X,Q}=\lambda(g_X)\,S_{X,Q}.
\end{equation}
Entries with $g_X=1$ and $S_{V_\ev,Q}$ do not change ($\lambda(1)=1$). 
The row module $X$ does not require normalization.

\begin{remark} 
\label{remconv}
 Let us make remarks concerning conventions and comparison with \cite{DL}.  
We distinguish $S_{X,Q}$ (row $X$), its complex conjugate $\ol{S_{X,Q}}$, 
and the entry $(S^{-1})_{Q,X}$ of the inverse of $\rho(S)$. 
 For the blocks used here $(S^{-1})_{Q,X}=\ol{S_{X,Q}}$ by Proposition \ref{propunit}.
 Symmetry is used for $S^U$ and the untwisted $S$ of $\Vev$ and 
 never for entries with a twisted index. 
 For $G$ of prime order $p$, Lemma \ref{lemorbrel} (b) reads 
$S^U_{M^1_{\Lambda_s},M^2_{\Lambda_t}}
 =\frac1pS_{M^1,M^2}\ol{\Lambda_s(g_2)}\Lambda_t(g_1^{-1})$. 
 Our conjugation, the argument $g_1^{-1}$, 
and the factor $1/p=|G_X|^{-1}$ coincide with those of \cite{DL} with $\Lambda_s=\psi_s$.  
 The Verlinde formula \cite{DL} also carries the conjugate on the third factor as ours. 

Let us explain why the conjugate matters. 
 Without it the summand of Theorem \ref{thmTV} would be multiplied by 
$\lambda(g_1g_2g_3)=\lambda(g_3^2)$ under a renormalization of the columns 
by \eqref{eqnormchange}.  
 In Example \ref{exA}, $n=3$, for $(T^{(1)}_0,T^{(1)}_0,T^{(2)}_0)$, 
 whose true fusion coefficient is one, the sum with the conjugate equals one, 
 and remains one after multiplying all column normalizations by 
 $\lambda(\sigma)=e^{2\pi i/3}$, whereas the sum without it equals 
 $1/4$ and becomes $-0.125+0.217\,i$ 
 (scripts \texttt{verify\_all.py}, \texttt{verify\_conv.py}). 

 In Example \ref{exA} the twisted trace functions are  theta functions. 
For $n=3$, $9$ and all sectors $a\ne0$, the transformation law \eqref{eqlaw} 
with the entries $S_{T^{(a)}_m,Q_k}=\frac12e^{-2\pi imk/4}e^{-2\pi iak/2n}$ 
of Lemma \ref{lemorbrel} (computed from $S^U$) holds to $10^{-15}$ 
at three values of $\tau$ 
(script \texttt{verify\_twisted\_traces.py}). 
 Reversing the sign of $a$ or taking complex conjugates produces errors of order 
 $0.7$-$1.2$. This tests the identification of the entries with the transformation 
 law of twisted traces, including the placement of $g$ versus $g^{-1}$ and 
 the normalization $\varphi_{Q_k}=e^{-2\pi ik/2n}\sigma$.
\end{remark}

\subsection{The twisted Verlinde formula for $\Vev$ for cyclic groups}
\label{secTV}

For $M_i\in\M(g_i)$, $g_3=g_1g_2$, set
\[
\mathsf T(M_1,M_2,M_3):=\sum_{Q\in\Q}\frac{S_{M_1,Q}\,S_{M_2,Q}\,\ol{S_{M_3,Q}}}{S_{V_\ev,Q}}.
\]
By \eqref{eqnormchange}, each summand is multiplied by 
$\lambda(g_1)\lambda(g_2)\ol{\lambda(g_3)}=\lambda(g_1g_2g_3^{-1})=1$ under a change of 
the normalizations of $Q$. This also holds if some $g_i=1$, where $\lambda(g_i)=1$.   
 Thus, $\mathsf T$ is well defined, and it is constant on orbits in each argument by Lemma 
\ref{lemorbinv} (a).
 The theorem below provides the first main result of this paper,
 the twisted Verlinde formula for $\Vev$. 

\begin{theorem} 
\label{thmTV}
Adopt both assumptions of Subsection \ref{secsetting}.
  Let $g_1$, $g_2\in G$ with $\langle g_1,g_2\rangle=G$, $g_3=g_1g_2$, 
 $M_i\in\M(g_i)$ with orbits $O_i$. Then:
\begin{enumerate}

\item[(i)] $\displaystyle \frac1{|O_1||O_2|}\sum_{N_1\in O_1,\,N_2\in O_2}N^{M_3}_{N_1,N_2}
=\mathsf T(M_1,M_2,M_3)$.

\item[(ii)] If at least two of $M_1$, $M_2$, $M_3$ are $G$-stable, 
 then $N^{M_3}_{M_1,M_2}=\mathsf T(M_1,M_2,M_3)$.

\end{enumerate}
\end{theorem}

Part (ii) always applies if $n$ is a prime power, and if every irreducible twisted 
$\Vev$-module is $G$-stable (e.g., for inner automorphisms,  Section \ref{secinner}). 
For $n$ prime it is given in \cite{DL}.

\begin{proof}
In this proof $\fus$ denotes the tensor product of $\Rep(\Vev)$ 
and $\fus_U$ that of $\mathcal C_U$. 
 No tensor product of twisted $V$-modules occurs (see Definition \ref{deffusion}).

\emph{0) Reduction to $U$ and multiplicities:} by Lemma \ref{lemF0} and Proposition 
\ref{factrep} (e), for arbitrary $\chi_1\in\widehat{G_{M_1}}$, $\chi_2\in\widehat{G_{M_2}}$, 
\[
\bigoplus_{N_1\in O_1,N_2\in O_2}N_1\fus N_2=F_0(M_{1,\chi_1})\fus F_0(M_{2,\chi_2})
=F_0(M_{1,\chi_1}\fus_UM_{2,\chi_2})=\bigoplus_YN^Y_{M_{1,\chi_1},M_{2,\chi_2}}F_0(Y),
\]
where $Y$ runs over the irreducible $U$-modules and each $F_0(Y)=\bigoplus_{N\in O_Z}N$, 
 $Y=Z_\psi$,  is multiplicity free (see Lemma \ref{lemF0}).
  Fix the simple object $M_3$. 
 It occurs in $F_0(Y)$ iff $M_3\in O_Z$, i.e., iff $Y\simeq(M_3)_{\chi_3}$ for some 
$\chi_3\in\widehat{G_{M_3}}$ (Proposition \ref{factclass} (iii)), and then exactly once.  
 The modules $(M_3)_{\chi_3}$, $\chi_3\in\widehat{G_{M_3}}$, are pairwise inequivalent 
(Proposition \ref{factclass} (i)). Thus, the $U$-modules $Y$ whose induction contains 
$M_3$ are these $|G_{M_3}|$ modules, each contributing once. 
 Comparing the multiplicities of $M_3$ we obtain 
\begin{equation}
\label{eqE1}
\sum_{N_1\in O_1,N_2\in O_2}N^{M_3}_{N_1,N_2}
=\sum_{\chi_3\in\widehat{G_{M_3}}}N^{(M_3)_{\chi_3}}_{M_{1,\chi_1},M_{2,\chi_2}}.
\end{equation}

\emph{1) Verlinde formula for $U$:}  
 by Huang version of Verlinde formula \cite{H1} for $\mathcal C_U$,
$N^{Y_3}_{Y_1,Y_2}=\sum_YS^U_{Y_1,Y}S^U_{Y_2,Y}\ol{S^U_{Y_3,Y}}/S^U_{U,Y}$. 
Apply it with $Y_i=(M_i)_{\chi_i}$. By Lemma \ref{lemorbrel} (a), 
a summand with $Y=N_\psi$ vanishes unless $g_1$, $g_2$, $g_3\in G_N$. 
 Since $G_N$ is a subgroup and $\langle g_1,g_2\rangle=G$ this leads to 
$G_N=G$, i.e., $Y=W_t$ with $W\in\M(h)$ stable, $h\in G$, $t\in\Z/n$.

\emph{2) Evaluation, characters and domains:} fix $\chi_1$, $\chi_2$. For stable $W\in\M(h)$ 
and $t\in\Z/n$, Lemma \ref{lemorbrel} (b) applies to the rows $(M_i)_{\chi_i}$, $i=1$, $2$,
 and to the complex conjugate of the row $(M_3)_{\chi_3}$, i.e.,  
$\ol{S^U_{(M_3)_{\chi_3},W_t}}
=|G_{M_3}|^{-1}\chi_3(h)\,\ol{\psi_t(g_3^{-1})}\;\ol{S^{[h]}_{M_3,W}}$.
  Moreover, $S^U_{U,W_t}=n^{-1}S^{[h]}_{V_\ev,W}$ for $X=\Vev$, $\chi=\psi_0$, $g=1$. 
 The $(W,t)$-summand therefore vanishes unless $h\in G_{M_1}\cap G_{M_2}\cap G_{M_3}$, and in that case equals
\begin{equation}\label{eqsummand}
\frac{\prod_{i=1}^3|G_{M_i}|^{-1}\;\ol{\chi_1(h)}\,\ol{\chi_2(h)}\,\chi_3(h)\;\;\psi_t(g_1^{-1})\psi_t(g_2^{-1})\ol{\psi_t(g_3^{-1})}\;\;S^{[h]}_{M_1,W}\,S^{[h]}_{M_2,W}\,\ol{S^{[h]}_{M_3,W}}}{n^{-1}\,S^{[h]}_{V_\ev,W}} .
\end{equation}

(i) \emph{ For $t$-sum:} the $t$-dependence is 
$\psi_t(g_1^{-1})\psi_t(g_2^{-1})\ol{\psi_t(g_3^{-1})}=\psi_t(g_1^{-1}g_2^{-1}g_3)
=\psi_t(1)=1$ since $g_3=g_1g_2$. 
  Hence, $\sum_{t\in\Z/n}$ gives $n$.

(ii) \emph{For $\chi_3$-sum:} the dependence on $\chi_3\in\widehat{G_{M_3}}$ 
 is in the factor $\chi_3(h)$ with $h\in G_{M_3}$. 
 We use the orthogonality of the cyclic group $G_{M_3}$ characters 
 of order $|G_{M_3}|$, i.e., 
 $\sum_{\chi_3\in\widehat{G_{M_3}}}\chi_3(h)=|G_{M_3}|$ 
 if $h=1$ and $0$ if $h\in G_{M_3}\setminus\{1\}$.  
 Terms with $h\notin G_{M_3}$ vanish already. 
 Thus, only $h=1$, i.e., $W\in\Q$, contributes. 
 Then $S^{[1]}=S$ and all character factors in \eqref{eqsummand} equal one. 

(iii) Summing \eqref{eqsummand} over $t$, $\chi_3$ and $W\in\Q$, we obtain 
\[
\sum_{\chi_3}N^{(M_3)_{\chi_3}}_{M_{1,\chi_1},M_{2,\chi_2}}=\frac{n\cdot n\cdot|G_{M_3}|}{|G_{M_1}||G_{M_2}||G_{M_3}|}\,\mathsf T(M_1,M_2,M_3)=|O_1||O_2|\,\mathsf T(M_1,M_2,M_3),
\]
since $n/|G_{M_i}|=|O_i|$ with one factor $n$ comming from $\sum_t$, 
the other from $S^U_{U,W_t}=n^{-1}S_{V_\ev,W}$). Together with \eqref{eqE1} this gives (i).

\emph{3) Individual coefficients:} if $M_1$, $M_2$ are stable, 
$|O_1|=|O_2|=1$, and (i) holds.  
 If $M_2$, $M_3$ are stable, then for $a\in G$, by Proposition \ref{factrep} (d) 
applied with $h=a^{-1}$, 
$N^{M_3}_{M_1\circ a,M_2}=N^{M_3\circ a^{-1}}_{M_1,M_2\circ a^{-1}}=N^{M_3}_{M_1,M_2}$ 
since $M_2\circ a^{-1}\simeq M_2$ and $M_3\circ a^{-1}\simeq M_3$. 
Thus, the left hand side of \eqref{eqE1} equals $|O_1|N^{M_3}_{M_1,M_2}$, with $|O_2|=1$. 
 The case when $M_1$, $M_3$ are stable is symmetric, i.e., 
  $N^{M_3}_{M_1,M_2\circ b}=N^{M_3\circ b^{-1}}_{M_1\circ b^{-1},M_2}=N^{M_3}_{M_1,M_2}$.
  Finally, let $n=p^a$. The residues of $g_1$, $g_2$, $g_3=g_1g_2$ modulo $p$ are 
  $(r_1,r_2,r_1+r_2)$ with $(r_1,r_2)\ne(0,0)$, as $\langle g_1,g_2\rangle=G$. 
 If $r_1=0$ then $g_2$, $g_3$ have residue $r_2\ne0$. 
 If $r_2=0$ then $g_1$, $g_3$ do.  
 If $r_1$, $r_2\ne0$ and $r_1+r_2=0$ then $g_1$, $g_2$ are generators. 
 In the remaining case all three are. 
 In each case at least two of $g_1$, $g_2$, $g_3$ generate $G=\Z/p^a$, 
 and a module twisted by a generator $g$ is stable since $g\in G_M$, \cite{DL}.
\end{proof}

\begin{remark}
\label{remn4}
Take $n=4$ and $(g_1,g_2,g_3)=(\sigma^2,\sigma,\sigma^3)$, 
with $M_1\in\M(\sigma^2)$ non-stable, 
 ($G_{M_1}=\langle\sigma^2\rangle$, $|O_1|=2$) and $M_2$, $M_3$ twisted by generators 
 and, hence, stable. 
 In 2), $h$ runs over $G_{M_1}\cap G_{M_2}\cap G_{M_3}=\langle\sigma^2\rangle$.
  $\chi_3$ runs over the four characters of $G_{M_3}=G$ and 
 $\sum_{\chi_3}\chi_3(h)=4\delta_{h,1}$. 
 The prefactor in (iii) is $\frac{n^2}{8}=2=|O_1||O_2|$. In 3), 
  the case when $M_2$, $M_3$ are stable, 
 the left side of \eqref{eqE1} is 
$N^{M_3}_{M_1,M_2}+N^{M_3}_{M_1\circ\sigma,M_2}=2N^{M_3}_{M_1,M_2}$. Hence,
 $N^{M_3}_{M_1,M_2}=\mathsf T(M_1,M_2,M_3)$. 
\end{remark}

\begin{remark}
\label{remTVopen}
 \emph{The remining open quesion:} for $n$ divisible by two distinct primes there are triples 
 with at most one stable module (e.g., $n=6$, $(g_1,g_2,g_3)=(\sigma^2,\sigma^3,\sigma^5)$, 
  $M_1$, $M_2$ are unstable). 
  There (i) gives only the orbit average over $O_1\times O_2$. 
  Whether the individual coefficients are given by $\mathsf T$ is not decided here 
 (by Lemma \ref{lemorbinv} (a) the summand $\mathsf T$ is constant on orbits,
  thus the claim would imply that $N^{M_3}_{M_1,M_2\circ b}$ does not depend on $b$).
  For non-cyclic abelian groups, projective representations of stabilizers appear 
  and the argument must be modified.
\end{remark}

Next proposition, the second main result of this paper,
 describes unitarity for generating sectors.  

\begin{proposition} 
\label{propunit}
Adopt both assumptions of Subsection \ref{secsetting} and let $g$ generate $G$. 
Then $|\M(g)|=|\Q|$ and the square matrix $(S_{X,Q})_{X\in\M(g),\,Q\in\Q}$ is unitary. 
In particular, for $X$, $X'\in\M(g)$ and $Q$, $Q'\in\Q$, 
\[
\sum_{Q\in\Q}S_{X,Q}\,\ol{S_{X',Q}}=\delta_{X,X'},
\qquad \sum_{X\in\M(g)}\ol{S_{X,Q}}\,S_{X,Q'}=\delta_{Q,Q'}.
\]
More generally $\sum_{N\in\M(h),\,\text{stable}}S^{[h]}_{X,N}\ol{S^{[h]}_{X',N}}=\delta_{X,X'}$
 for every $h\in G$. 
\end{proposition}

\begin{proof}
\emph{Orthonormality of rows:} let $X$, $X'\in\M(g)$. They are stable,  
 $g \in G_X$, $\langle g\rangle=G$.  
 Thus, all $n$ components $X_s$, $X'_{s'}$ exist. Unitarity of $S^U$ gives 
 $\sum_YS^U_{X_s,Y}\ol{S^U_{X'_{s'},Y}}=\delta_{XX'}\delta_{ss'}$ since  
  distinct stable modules in $\M(g)$ give inequivalent $U$-modules, 
 see Proposition \ref{factclass} (iii).  
 By Lemma \ref{lemorbrel} (a) only $Y$ with $g\in G_Y$, i.e., $G_Y=G$, contribute: 
 $Y=N_t$, $N\in\M(h)$ stable, $h\in G$, $t\in\Z/n$, each such $Y$ once. 
 By Lemma \ref{lemorbrel} (b) with $G_X=G_{X'}=G$, 
 $S^U_{X_s,N_t}=n^{-1}\ol{\psi_s(h)}\psi_t(g^{-1})S^{[h]}_{X,N}$. 
 Thus, the left hand side equals
\[
n^{-2}\sum_{h\in G}\ol{\psi_s(h)}\psi_{s'}(h)\,A_h\sum_{t\in\Z/n}|\psi_t(g^{-1})|^2
=n^{-1}\sum_{h\in G}\psi_{s'-s}(h)A_h,
\qquad A_h:=\sum_{N\in\M(h)\text{ stable}}S^{[h]}_{X,N}\ol{S^{[h]}_{X',N}} ,
\]
where complex conjugation acts on the second factor of each product. 
For $X=X'$ the left side is $\delta_{ss'}$ and for $X\ne X'$ it is zero. 
 With $a=s'-s$ ranging over $\Z/n$ this gives 
 $\sum_h\psi_a(h)A_h=n\delta_{XX'}\delta_{a,0}$ for all $a$, 
 hence, $A_h=\delta_{XX'}$ for all $h\in G$ by Fourier inversion on $G$, 
 i.e., $A_h=n^{-1}\sum_a\ol{\psi_a(h)}\,n\delta_{XX'}\delta_{a,0}=\delta_{XX'}$. 
 For $h=1$, $A_1=\sum_{Q\in\Q}S_{X,Q}\ol{S_{X',Q}}$.

 Denote by $W_{(g,h)}$ the span of the functions $Z_M(v,(g,h),\tau)$, $M\in\M(g,h)$
  with $v\in\Vev$, as in \eqref{eqlaw}. By \cite{DLM3,DL}, the functions of all sectors 
 together are linearly independent. The dimension of their span is $\sum_{(g,h)}|\M(g,h)|$.   
 In particular, $\dim W_{(g,h)}=|\M(g,h)|$ for each commuting pair, with $\M(g,h)$ 
 the set of isomorphism classes of individual modules, 
 so that, for sectors where some module is not $h$-stable or has non-trivial orbit, 
 the count is that of \cite{DL}.
 By the transformation law, $\rho(S)$ maps $W_{(g,h)}$ into $W_{(h,g^{-1})}$, and 
 $\rho(S)^{-1}=\rho(S^{-1})$ maps $W_{(h,g^{-1})}$ into $W_{(g,h)}$.  
 Indeed, we apply the law to 
 $S^{-1}=\left(\begin{smallmatrix}0&1\\-1&0\end{smallmatrix}\right)$, 
 which sends $(h,g^{-1})$ to $(g,h)$. Thus, $\rho(S)$ restricts to a linear 
 isomorphism $W_{(g,h)}\to W_{(h,g^{-1})}$, and $|\M(g,h)|=|\M(h,g^{-1})|$. 
 For $h=1$ and $g$ a generator this reads $|\M(g)|=|\M(1,g)|$. 
 Here $\M(g,1)=\M(g)$ (all modules are stable, thus orbits are singletons) 
 and $\M(1,g)=\Q$ since $g$ generates $G$.  
 A square matrix with orthonormal rows is unitary, thus its columns are orthonormal.
\end{proof}


Let $A=\Vev\oplus J$ be the superalgebra in $\Rep(\Vev)$. 
For a simple $M\in\Rep(\Vev)$ set $F_V(M)=M\oplus J\fus M$ with $M$ even and $J\fus M$ odd. 
If $\{X,J\fus X\}$ is an orbit we call $X^\ev:=X$ an \emph{even representative}. 
 Then $F_V(J\fus X)=\Pi F_V(X)$, the same module with opposite grading. 
 In this paper, $N^{X_3}_{X_1,X_2}$ always denotes a fusion coefficient of $\Rep(\Vev)$ 
 (see Definition \ref{deffusion}). Quantities for $V$-modules appear under the 
 assumptions below.
  
\subsubsection{Minimal assumptions (Theorem \ref{thmmain} (b)):}
\label{hypdict} 

\emph{}

  \emph{First}, for every simple $M\in\Rep(\Vev)$ of twist $g$, $F_V(M)$ is a 
$\mathbb Z/2\mathbb Z$-graded $V$-module, $g$-twisted if $\epsilon_M=1$ and 
$gP$-twisted if $\epsilon_M=-1$ (see Definition \ref{defsigns}), 
 this includes $J$-fixed $M$. 

 \emph{Second}, For simple $X_1$, $X_2$, $X_3$ with $g_3=g_1g_2$ and even representatives 
$X_i^\ev=X_i$, the space $I(F_V(X_3);F_V(X_1),F_V(X_2))$ of intertwining operators 
is $\mathbb Z/2\mathbb Z$-graded (an operator $\mathcal Y$ has parity $a$ if 
$\mathcal Y(u,z)w\in F_V(X_3)^{|u|+|w|+a}$), 
and its even part has dimension $N^{X_3}_{X_1,X_2}$, its odd part $N^{J\fus X_3}_{X_1,X_2}$. 

\subsubsection{Classification and monoidality (Corollary \ref{corproduct}:}
\label{hypmono}

\emph{}

\emph{First:} if $J\fus M\not\simeq M$, $F_V(M)$ is simple as a graded module, and every simple 
graded $g$- or $gP$-twisted $V$-module $N$ with $N^{\bar 0}\not\simeq J\fus N^{\bar 0}$ 
is isomorphic to some $\Pi^iF_V(M)$ for a unique free orbit.

\emph{Second:} on the full subcategory $\Rep(V)^{\rm free}$ generated by the $\Pi^iF_V(M)$, $M$ free, 
a tensor product $\fus_V$ exists, and $F_V(M)\fus_VF_V(M')\simeq F_V(M\fus M')$ whenever every
 simple part of $M\fus M'$ is $J$-free.

\subsubsection{Li transport for superalgebras (Section \ref{secinner}:)} 
\label{hypL} 

Let $h\in\Vev$ be as in  Section \ref{secinner}.

\emph{First:}  
 for every irreducible (parity-untwisted, parity-twisted, or twisted) $V$-module $M$, 
 $M^{(h)}=(M,Y_M(\Delta(h,z)\cdot,z))$ is a module twisted by $\phi(h)=e^{-2\pi ih_0}$ 
 (composed with the same parity twist), and $\Delta(h_1,z)\Delta(h_2,z)=\Delta(h_1+h_2,z)$.

\emph{Second:}  
an intertwining operator $\mathcal Y$ of type $\binom{M_3}{M_1\,M_2}$ gives an 
intertwining operator $\mathcal Y(\Delta(h,z)\cdot,z)$ of type 
$\binom{M_3^{(h)}}{M_1\ M_2^{(h)}}$.

\emph{Third:} 
 (Two twisted arguments.) $N^{M_3^{((a+b)h)}}_{M_1^{(ah)},M_2^{(bh)}}=N^{M_3}_{M_1,M_2}$, 
where $M^{(nh)}$ denotes the module obtained by the full spectral flow, an untwisted module.
 The corresponding identification of spaces of intertwining operators respects the 
$\Z/2\Z$-gradings.

Note that conditions \ref{hypL}, \ref{hypmono}, \ref{hypdict}   
 are additional conditions are not consequences of two assumptions  
of Subsection \ref{secsetting} or of the cited results. 
 For $g=1$, conditions \ref{hypdict}-\ref{hypmono}  
 are established in \cite{HKL, CKM}.  
For vertex operator algebras, the first and second conditions of \ref{hypL} 
 are theorems \cite{Li, DLM1}.
  Assumptions \ref{hypL} state them for superalgebras.  
 The third condition of \ref{hypL} is not covered by those results 
(where only two of the three arguments are twisted).    
   For inner automorphisms,  
conditions \ref{hypdict}-\ref{hypmono}   
reduce to the untwisted theory under assumptions \ref{hypL}  
 (see Proposition \ref{propinner}). 

\section{Sign data and equivariance of the orbifold $S$-matrix}\label{secsign}

First, let us list conventions for this section: 

(i) For $W\in\{\Vev,U\}$, $S^W$ is the $S$-matrix of $\mathcal C_W$ defined by modular 
transformation (for $W=U$: of the functions $Z_Y(v,\tau)$, $v\in U$). An $S$ without  
superscript is the untwisted $S$-matrix of $\Vev$ or an orbifold entry $S_{X,Q}$ (known from 
the indices). 

(ii) $\theta_A=e^{2\pi ih_A}$. 

(iii) For simple $Z$, $Y$ in a modular category  
$\mathcal C_W$ set $\chi_Y(Z):=S_{Z,Y}/S_{\mathbf 0,Y}$. By the Verlinde formula this is a  
character of the fusion ring, thus $S^W_{A\fus B,Y}=\chi_Y(A)S^W_{B,Y}$. 

(iv) Symmetry is used only for $S^U$ and for the untwisted $S$ of $\Vev$. 

(v) The quantities $d$, $\kappa$, $\epsilon$, $\eta$, $\pi$ of Definition \ref{defsigns} 
are ratios of $S$-matrix entries. Their interpretation in terms of $V$-modules requires 
 conditions \ref{hypdict}.

The lemma below provides the Gauss identity. 

\begin{lemma}
\label{lemgauss}
Let $W$ be a simple, rational, $C_2$-cofinite, self-dual vertex operator algebra of CFT type 
(we use $W=U$ and $W=\Vev$) with irreducible modules $a$, $b$, $\dots$, vacuum module 
$\mathbf 0=W$, conformal weights $h_a$, twists $\theta_a=e^{2\pi ih_a}$, $S=S^W$ and fusion 
coefficients $N_{ab}^f$. Then  
\[
S_{ab}=\theta_a\,\theta_b\sum_fN_{ab}^f\,\theta_f^{-1}\,S_{f\mathbf 0}\qquad\text{for all }a, b.
\]
With the opposite convention for $T$ one replaces $\theta$ by $\theta^{-1}$, and 
 the central charge drops out.  
\end{lemma}

\begin{proof}

Let $\hat T=\mathrm{diag}(\gamma\theta_a)$, $\gamma=e^{-2\pi ic/24}$. 
 By Zhu theorem the characters of the irreducible $W$-modules are linearly independent and 
span an $\mathrm{SL}_2(\Z)$-representation in which $S$, $\hat T$ are the matrices of 
$\tau\mapsto-1/\tau$, $\tau\mapsto\tau+1$. Hence, $(S\hat T)^3=S^2$ or $(\hat TS)^3=S^2$, 
according to the order of composition. In the first case $S\hat TS\hat TS\hat T=SS$, 
and cancelling $S$ on the left gives $\hat TS\hat TS\hat T=S$. In the second, 
$\hat TS\hat TS\hat TS=SS$, and cancelling $S$ on the right gives the same relation. 
Multiplying $\hat TS\hat TS\hat T=S$ on the left and on the right by $\hat T^{-1}$ gives 
$S\hat TS=\hat T^{-1}S\hat T^{-1}$. Evaluate the $(a,b)$ entry: 
since $\hat T_e=\gamma\theta_e$, we have 
\begin{equation}\label{eqSTS}
\sum_eS_{ae}S_{be}\,\theta_e=\gamma^{-3}\,\theta_a^{-1}\theta_b^{-1}\,S_{ab}.
\end{equation}
By the Verlinde formula, $e\mapsto S_{ae}/S_{\mathbf 0e}$ is a character of the fusion ring, 
 thus $S_{ae}S_{be}/S_{\mathbf 0e}=\sum_fN_{ab}^fS_{fe}$. Hence, using \eqref{eqSTS} for 
the pair $(f,\mathbf 0)$ and $\theta_{\mathbf 0}=1$, we get 
\[
\sum_eS_{ae}S_{be}\theta_e=\sum_fN_{ab}^f\sum_eS_{fe}S_{\mathbf 0e}\theta_e
=\gamma^{-3}\sum_fN_{ab}^f\,\theta_f^{-1}S_{f\mathbf 0}.
\]
Comparing with \eqref{eqSTS} (the factor $\gamma^{-3}$ cancels) gives the statement of Lemma. 
\end{proof}

\begin{lemma}
\label{lemnonvan}
$S_{X,V_\ev}\ne0$ for every $X\in\R$, and $S_{V_\ev,Q}\ne0$ for every $Q\in\Q$.
\end{lemma}

\begin{proof}
In the modular category $\mathcal C_U$, $S^U_{U,Y}=S^U_{U,U}\dim_qY$ with 
$S^U_{U,U}\ne 0$ and $\dim_qY\ne 0$ for every simple $Y$, 
 since dimensions of simple objects of a fusion category are nonzero,  
and $S^U$ is symmetric.  Thus, $S^U_{Y,U}=S^U_{U,Y}\ne 0$. 
The column $V_\ev\in\Q$ with $t=0$ is $U$ (with normalization $\varphi_{V_\ev}=\sigma$). 
Lemma \ref{lemorbrel} (b) with $W=V_\ev$, $h=1$, $t=0$ gives, for every $X\in\M(g)$ and 
any $\chi$, 
$S^U_{X_\chi,U}=|G_X|^{-1}\ol{\chi(1)}\psi_0(g^{-1})S_{X,V_\ev}=|G_X|^{-1}S_{X,V_\ev}$, 
hence, $S_{X,V_\ev}\ne0$. With $X=\Vev$ ($g=1$, $\chi=\psi_0$, $X_\chi=U$) and column $Q$: 
$S^U_{U,Q_t}=n^{-1}S_{V_\ev,Q}$, which is nonzero.
\end{proof}

\begin{definition}
\label{defsigns}
With the normalizations $\varphi_J=\sigma|_J$, $\varphi_{V_\ev}=\sigma|_{\Vev}$ on the columns 
$J$, $V_\ev\in\Q$, set  
\[
d:=\frac{S_{J,V_\ev}}{S_{V_\ev,V_\ev}},\qquad \kappa_X:
=\frac{S_{X,J}}{S_{X,V_\ev}}\ (X\in\R),\qquad 
\epsilon_X:=\frac{\kappa_X}{d},\qquad \eta_Q:=\kappa_Q\ (Q\in\Q),
\]
and the parity type $\pi(X)\in\Z/2\Z$ by $(-1)^{\pi(X)}=\epsilon_X$. 
 All equalities and sums of values of $\pi$ are taken in $\Z/2\Z$.  
Thus, $\pi(X_1)+\pi(X_2)=\pi(X_3)$ iff $\epsilon_{X_1}\epsilon_{X_2}\epsilon_{X_3}=1$. 
These are ratios of $S$-matrix entries (well defined by Lemma \ref{lemnonvan}). 
Lemma \ref{lemvals} (2) identifies them with values of characters of $\mathcal C_U$, 
and conditions \ref{hypdict} gives $\epsilon_X$ a module-theoretic meaning (parity type). 
\end{definition}

Induced normalization and uniform shift are given in 
\begin{lemma}
\label{lemshift}
Let $X\in\M(g)$ with stabilizer $G_X=\langle\sigma_X\rangle$ of order $m$ and fixed $\varphi_X$. Then $J\fus X\in\M(g)$ has $G_{J\fus X}=G_X$, and:

\begin{enumerate}

\item[(i)] Let $0\ne\mathcal Y\in I(J\fus X;J,X)$ (a one-dimensional space).  
There is a unique isomorphism $\theta:J\fus X\to(J\fus X)\circ\sigma_X$ with
\begin{equation}
\label{eqind}
\theta\,\mathcal Y(u,z)w=\mathcal Y(\sigma_Xu,z)\,\varphi_X(\sigma_X)w, \qquad u\in J,\ w\in X,
\end{equation}
and $\theta^m=1$. Hence, $\varphi^{\rm ind}(\sigma_X^k):=\theta^k$ is a representation of 
$G_X$ on $J\fus X$ (the induced normalization), satisfying \eqref{eqind} for all 
$h\in G_X$.

\item[(ii)] With respect to $\varphi_X$ and $\varphi^{\rm ind}$: $J_0\fus_UX_\chi\simeq(J\fus X)_\chi$ \emph{for every} $\chi\in\widehat{G_X}$ (the label is not shifted).

\item[(iii)] If $J\fus X$ is given another normalization $\lambda\varphi^{\rm ind}$, $\lambda\in\widehat{G_X}$, then $J_0\fus_UX_\chi\simeq(J\fus X)'_{\chi\lambda}$ for all $\chi$: the shift is the same character $\lambda$ for all $\chi$.
\end{enumerate}
\end{lemma}

\begin{proof}
\emph{Stabilizer:} by Proposition \ref{factrep} (d), 
$(J\fus X)\circ a\simeq(J\circ a)\fus(X\circ a)\simeq J\fus(X\circ a)$ as $J$ is stable, 
and $J\fus(\cdot)$ is an equivalence. Thus, $J\fus(X\circ a)\simeq J\fus X$ iff 
$X\circ a\simeq X$.

\emph{One-dimensionality:} $J$ is invertible in $\Rep(\Vev)$, thus $J\fus X$ is simple and 
$N^{J\fus X}_{J,X}=1$. By Proposition \ref{factrep} (b), $\dim I(J\fus X;J,X)=1$. 

\emph{(i) Existence and uniqueness:} let $\Theta$ be the set of 
isomorphisms $\theta_0:J\fus X\to(J\fus X)\circ\sigma_X$. It is nonempty since $J\fus X$ 
is $G_X$-stable, and is a $\C^\times$-torsor. The same map $\mathcal Y$ is an intertwining 
operator of type $\binom{(J\fus X)\circ\sigma_X}{J\circ\sigma_X\ \ X\circ\sigma_X}$,  
 \cite{DL}. Explicitly, $M\circ h$ has 
vertex operators $Y_{M\circ h}(v,z)=Y_M(hv,z)$. The Jacobi identity for  
$\mathcal Y$ is, for each $v\in\Vev$ homogeneous for the automorphisms occurring in it 
(a joint eigenvector), an identity between $Y_{M_3}(v,\cdot)$, $Y_{M_2}(v,\cdot)$, 
$Y_{M_1}(v,\cdot)$ and $\mathcal Y$. Since $G$ is abelian, $hv$ belongs to the same 
joint eigenspace as $v$. 
  Thus, the identity for $hv$, valid for the original modules, 
is the identity for $v$ in the $h$-conjugated modules. 
  Composing with the isomorphisms 
$\varphi_J(\sigma_X)=\sigma_X|_J:J\to J\circ\sigma_X$ and 
$\varphi_X(\sigma_X):X\to X\circ\sigma_X$, and then composing with $\theta_0^{-1}$ gives
$\mathcal Y^{\theta_0}(u,z)w:=\theta_0^{-1}\mathcal Y(\sigma_Xu,z)\varphi_X(\sigma_X)w
\in I(J\fus X;J,X)$ which is nonzero. 
 Hence, $\mathcal Y^{\theta_0}=c(\theta_0)\mathcal Y$ with 
$c(\theta_0)\in\C^\times$. 
  For $\mu\in\C^\times$, $c(\mu\theta_0)=\mu^{-1}c(\theta_0)$. 
Thus $\theta:=c(\theta_0)\theta_0$ is the unique element of $\Theta$ with $c(\theta)=1$, 
which is \eqref{eqind}. Uniqueness can also be seen: two solutions differ by a scalar 
$\mu$ with $(\mu-1)\mathcal Y=0$. 

\emph{The relation $\theta^m=1$:} iterate \eqref{eqind}  
 by applying it to $(\sigma_Xu,\varphi_X(\sigma_X)w)$. The functor $\circ\sigma_X$ does 
not change underlying linear maps. Then it gives   
 $\theta^k\mathcal Y(u,z)w=\mathcal Y(\sigma_X^ku,z)\varphi_X(\sigma_X)^kw$ for $k\ge1$, 
where $\theta^k$ is the composite of $k$ copies of $\theta$, a linear isomorphism 
$J\fus X\to(J\fus X)\circ\sigma_X^k$ compatible with the module structures. For $k=m$ we 
have $\sigma_X^m=1$. Thus, $\theta^m$ is an endomorphism of the simple $\Vev$-module 
$J\fus X$, hence a scalar $\mu$ by Schur lemma. Since $\sigma_X^m=1$ on $J$ and 
$\varphi_X(\sigma_X)^m=\varphi_X(1)=1$ ($\varphi_X$ is a representation of $G_X$), 
$\theta^m$ fixes coefficients of $\mathcal Y(u,z)w$ which are not all zero, 
 because $\mathcal Y(u,z)w\ne 0$ for $u\ne0\ne w$ by \cite{DL}. 
Hence, $\mu=1$ and $\theta^m=\id$. 

\emph{(ii)} For $h\in G_X$, $u\in J_0$ ($hu=u$) and $w\in X_\chi$, \eqref{eqind} gives 
$\varphi^{\rm ind}(h)\mathcal Y(u,z)w=\chi(h)\mathcal Y(u,z)w$. Thus, 
$\mathcal Y|_{J_0\otimes X_\chi}$ is a nonzero $U$-intertwining operator with coefficients 
in $(J\fus X)_\chi$, an irreducible $U$-module (see Proposition \ref{factclass} (i)). 
Thus, $N^{(J\fus X)_\chi}_{J_0,X_\chi}\ge 1$. Since $J_0$ is a simple current of  
$\mathcal C_U$, $J_0\fus_UX_\chi$ is simple and $J_0\fus_UX_\chi \simeq(J\fus X)_\chi$.

\emph{(iii)} Multiplying $\varphi^{\rm ind}(h)$ by $\lambda(h)$ multiplies the eigenvalue 
on every $(J\fus X)_\chi$ by $\lambda(h)$, i.e., relabels simultaneously 
$\chi\mapsto\chi\lambda$ for all $\chi$.
\end{proof}

\begin{lemma} 
\label{lemvals}
\begin{enumerate}

\item $d=\dim_qJ=\dim_qJ_0\in\{\pm1\}$.

\item $\kappa_X=\chi_{X_\chi}(J_0)\in\{\pm1\}$ for every $X\in\R$ and every 
$\chi\in\widehat{G_X}$. Hence $\epsilon_X$, $\eta_Q\in\{\pm1\}$.

\item $\chi_{J_0}(J_0)=d$.

\item $\kappa_{J\fus X}=\kappa_X$. Hence, $\epsilon_{J\fus X}=\epsilon_X$, 
$\pi(J\fus X)=\pi(X)$, $\eta_{J\fus Q}=\eta_Q$. 

\end{enumerate}
\end{lemma}

\begin{proof}
(1) Lemma \ref{lemorbrel} (b) with rows $J$ ($g=1$, $\chi=\psi_0$, thus $X_\chi=J_0$) 
 and column $V_\ev$ ($t=0$): $S^U_{J_0,U}=n^{-1}S_{J,V_\ev}$, 
$S^U_{U,U}=n^{-1}S_{V_\ev,V_\ev}$.  
Thus, $d=S^U_{J_0,U}/S^U_{U,U}=\dim_qJ_0$. Also, $d=S_{J,V_\ev}/S_{V_\ev,V_\ev}=\dim_qJ$ 
in $\mathcal C_{\Vev}$. As $J_0\fus_UJ_0\simeq U$, $d^2=1$.

(2) Lemma \ref{lemorbrel} (b) with columns $J$ and $V_\ev$ at $t=0$ 
(where $J_0$ and $U$ are the $t=0$ components) and the symmetry of $S^U$:
$\kappa_X=\dfrac{S^U_{X_\chi,J_0}}{S^U_{X_\chi,U}}=\dfrac{S^U_{J_0,X_\chi}}{S^U_{U,X_\chi}}
=\chi_{X_\chi}(J_0)$. Then $\chi(J_0)^2=\chi(J_0\fus_UJ_0)=\chi(U)=1$.

(3) Apply Lemma \ref{lemgauss} to $W=U$, $a=b=J_0$: since $J_0\fus_UJ_0\simeq U$ 
and $\theta_{J_0}^2=e^{4\pi ih_J}=1$ ($2h_J\in\Z$ because $V$ is $\tfrac12\Z$-graded), 
$S^U_{J_0,J_0}=S^U_{U,U}$. By (1), $S^U_{U,J_0}=S^U_{J_0,U}=d\,S^U_{U,U}$. 
Hence, $\chi_{J_0}(J_0)=S^U_{J_0,J_0}/S^U_{U,J_0}=1/d=d$.

(4) By Lemma \ref{lemshift} (ii), $(J\fus X)_\chi=J_0\fus_UX_\chi$ (with induced normalization). 
 By (2) $\kappa_{J\fus X}=\chi_{J_0\fus_UX_\chi}(J_0)
=S^U_{J_0,J_0\fus_UX_\chi}/S^U_{U,J_0\fus_UX_\chi}$. By symmetry, multiplicativity and (3):
$S^U_{J_0,J_0\fus X_\chi}=S^U_{J_0\fus X_\chi,J_0}
=\chi_{J_0}(J_0)S^U_{X_\chi,J_0}=d\,S^U_{J_0,X_\chi}$ and $S^U_{U,J_0\fus X_\chi}
=S^U_{J_0\fus X_\chi,U}=\chi_U(J_0)S^U_{X_\chi,U}=d\,S^U_{U,X_\chi}$. 
Dividing gives $\kappa_X$.
\end{proof}

Equivariance in the row index is given in 

\begin{lemma}
\label{lemrow}

For all $X\in\R$ and $Q\in\Q$: $\;S_{J\fus X,Q}=\eta_Q\,S_{X,Q}$.

\end{lemma}

\begin{proof}

Let $g=g_X$; $J\fus X$ has the same twist and stabilizer. Use the induced normalization on 
$J\fus X$, the row entries do not depend on it. By Lemma \ref{lemshift} (ii), 
$(J\fus X)_\chi=J_0\fus_UX_\chi$, thus 
$S^U_{(J\fus X)_\chi,Q_t}=\chi_{Q_t}(J_0)\,S^U_{X_\chi,Q_t}=\eta_Q\,S^U_{X_\chi,Q_t}$ 
(Lemma \ref{lemvals} (2) for $Q$). Lemma \ref{lemorbrel} (b) on both sides 
(column $Q$ untwisted and the factors $|G_X|^{-1}\psi_t(g^{-1})$ are the same) gives 
$S_{J\fus X,Q}=\eta_QS_{X,Q}$. There is no distinction between stable and unstable $X$  
since the stabilizer enters only through $|G_X|$. 
\end{proof}

 Equivariance in the column index is given by 

\begin{lemma}
\label{lemcol}

Let $Q\in\Q$, and give $J\fus Q\in\Q$ the induced normalization. Then for every $X\in\R$ 
(any stabilizer)
\[
S_{X,J\fus Q}=\kappa_X\,S_{X,Q},\qquad\text{in particular, }\; S_{V_\ev,J\fus Q}=d\,S_{V_\ev,Q}.
\]
If $J\fus Q$ carries the normalization $\lambda\varphi^{\rm ind}$, $\lambda\in\widehat G$, 
then $S_{X,J\fus Q}=\lambda(g_X)\kappa_XS_{X,Q}$.
\end{lemma}

\begin{proof}

Symmetry is applied here only to $S^U$. By Lemma \ref{lemshift} (ii), 
$(J\fus Q)_t=J_0\fus_UQ_t$ for all $t$. For $X\in\M(g)$ and $\chi\in\widehat{G_X}$,
\[
S^U_{X_\chi,(J\fus Q)_t}=S^U_{J_0\fus_UQ_t,X_\chi}
=\chi_{X_\chi}(J_0)\,S^U_{Q_t,X_\chi}=\kappa_X\,S^U_{X_\chi,Q_t}
\]
(by Lemma \ref{lemvals} (2)). Lemma \ref{lemorbrel} (b) with the column $J\fus Q$, 
respectively, $Q$, 
gives $|G_X|^{-1}\psi_t(g^{-1})S_{X,J\fus Q}=\kappa_X|G_X|^{-1}\psi_t(g^{-1})S_{X,Q}$. For a
 different normalization $\lambda\varphi^{\rm ind}$ of $J\fus Q$, \eqref{eqnormchange} 
multiplies the left side by $\lambda(g_X)$. Here $\lambda$ depends only on the normalization 
of the column $J\fus Q$ and does not depend on $X$. For three rows with twists 
$g_1$, $g_2$, $g_3=g_1g_2$ the 
resulting phases $\lambda(g_1)\lambda(g_2)\ol{\lambda(g_3)}=\lambda(g_1g_2g_3^{-1})=1$ cancel 
in $f_0$ (Theorem \ref{thmmain}, Part 2).

\end{proof}

For $g=1$ these are \cite{CKM} for $\mathcal C_{\Vev}$. For $g\ne1$, once normalizations are 
compatible no phase appears. With an arbitrary normalization the phase is the character 
$\lambda(g_X)$, which is the same for all $X$ and cancels in the fusion rules since 
$g_3=g_1g_2$ (see proof of Theorem \ref{thmmain}).

\begin{lemma}
\label{lemfree}
If $d=-1$ then $J\fus X\not\simeq X$ for every $X\in\R$.
\end{lemma}

\begin{proof}
If $J\fus X\simeq X$, Lemma \ref{lemrow} with $Q=V_\ev$ ($\eta_{V_\ev}=\kappa_{V_\ev}=d$) 
gives $S_{X,V_\ev}=dS_{X,V_\ev}=-S_{X,V_\ev}$, contradicting Lemma \ref{lemnonvan}.
\end{proof}

Concerning fixed points we state the following 
\begin{lemma}
\label{lemfixed}

\begin{enumerate}
\item[(i)] If $Q\in\Q$ and $J\fus Q\simeq Q$, then $d=1$ and $\eta_Q=-1$.

\item[(ii)] If $X\in\R$ and $J\fus X\simeq X$, then $S_{X,Q}=0$ for every $Q\in\Q$ with 
$\eta_Q=-1$.

\item[(iii)] If $g_X$ generates $G$, then $J\fus X\simeq X$ if and only if $S_{X,Q}=0$ for all 
$Q\in\Q$ with $\eta_Q=-1$, i.e., iff $\sum_Q\eta_Q|S_{X,Q}|^2=1$.

\end{enumerate}
\end{lemma}

\begin{proof}
(i) $d=1$ by Lemma \ref{lemfree}, thus $\theta_J=-1$ (Proposition \ref{factJ} (ii)). 
Apply Lemma \ref{lemgauss} to $W=\Vev$, $a=J$, $b=Q$ (an untwisted module). Then 
  $N_{JQ}^f=\delta_{f,J\fus Q}=\delta_{f,Q}$, 
thus $S_{J,Q}=\theta_J\theta_Q\theta_Q^{-1}S_{Q,V_\ev}=\theta_J\,S_{Q,V_\ev}$. 
By symmetry of the untwisted $S$ of $\Vev$, $\eta_Q=S_{Q,J}/S_{Q,V_\ev}=\theta_J=-1$.

(ii) By Lemma \ref{lemrow}, $S_{X,Q}=S_{J\fus X,Q}=\eta_QS_{X,Q}$. 

(iii) $X$ and $J\fus X$ both belong to $\M(g_X)$ (the same twist), 
 and $g_X$ generates $G$. Proposition \ref{propunit} (for rows) gives 
\[
\sum_{Q\in\Q}S_{J\fus X,Q}\ol{S_{X,Q}}=\delta_{J\fus X,X}, 
 \qquad\sum_{Q\in\Q}|S_{X,Q}|^2=1 .
\]
By Lemma \ref{lemrow}, $S_{J\fus X,Q}=\eta_QS_{X,Q}$, thus the first sum equals
 to $\sum_Q\eta_Q|S_{X,Q}|^2$. If $J\fus X\simeq X$, both sums equal one. 
 By subtracting, $\sum_Q(1-\eta_Q)|S_{X,Q}|^2=0$, i.e., 
$2\sum_{\eta_Q=-1}|S_{X,Q}|^2=0$. 
 Thus, $S_{X,Q}=0$ whenever $\eta_Q=-1$. 
 Conversely, if $S_{X,Q}=0$ whenever $\eta_Q=-1$ then 
$\sum_Q\eta_Q|S_{X,Q}|^2=\sum_Q|S_{X,Q}|^2=1$, 
thus $\delta_{J\fus X,X}=1$ and $J\fus X\simeq X$. 
  Note that (i) concerns untwisted fixed modules ($Q\in\Q$), and (iii) concerns 
twisted and untwisted $X$ with generating twist. 
\end{proof}

Conformal weight criterion in untwisted sector is given in 

\begin{lemma} 
\label{lemweight}

For every untwisted $X\in\M(1)$ (stable or not) 
$\theta_{J\fus X}/\theta_X=-d\,\epsilon_X$, i.e., 
$e^{2\pi i(h_{J\fus X}-h_X)}=-d\,\epsilon_X$ for the conformal weights.
\end{lemma}

\begin{proof}

Apply Lemma \ref{lemgauss} to $W=\Vev$, $a=J$, $b=X$: 
$N_{JX}^f=\delta_{f,J\fus X}$, thus 
$S_{J,X}=\theta_J\theta_X\theta_{J\fus X}^{-1}S_{J\fus X,V_\ev}$. 
Here $S_{J\fus X,V_\ev}=\chi_{V_\ev}(J)S_{X,V_\ev}=d\,S_{X,V_\ev}$ 
(multiplicativity of the character $\chi_{V_\ev}$ and $d=S_{J,V_\ev}/S_{V_\ev,V_\ev}$). 
  By symmetry of the untwisted $S$, $\kappa_X=S_{X,J}/S_{X,V_\ev}=S_{J,X}/S_{V_\ev,X}$, 
and, therefore, $\kappa_X=d\,\theta_J\,\theta_X/\theta_{J\fus X}$, i.e., 
$\theta_{J\fus X}/\theta_X=d\,\theta_J/\kappa_X=\theta_J\,\epsilon_X=-d\,\epsilon_X$ 
(using $d^2=\epsilon_X^2=1$, $\theta_J=-d$). All entries $S_{X,V_\ev}$ have cancelled, 
thus the result does not depend on the normalization of $S$. With the opposite convention 
for $T$ the computation gives the same result since $\epsilon_X\theta_J=\pm1$ is real.
\end{proof}

\begin{remark}
\label{remweights}
Let us make a remark about parity type and weights for fixed twisted modules. 
For $d=1$ ($J$ has weights in $\Z+\frac12$) Lemma \ref{lemweight} gives $\epsilon_X=1$ 
iff the weights of $J\fus X$ and $X$ differ by $\frac12$ mod $\Z$ (NS-type); 
 $\epsilon_X=-1$ iff they agree mod $\Z$ (Ramond-type). 
 For $d=-1$ the roles are reversed. 
 For twisted $X$ and $n$ odd, the analogous statement 
 $e^{2\pi in(h_{J\fus X}-h_X)}=-d\,\epsilon_X$ holds in Examples \ref{exA} and \ref{exC}.  
    The first condition of \ref{hypdict} requires  that for the sign $\epsilon_X$,  
 and it is an independent check of the labelling $\epsilon\leftrightarrow$ 
 parity type that does not use $S$-matrices. 
    Note that an untwisted fixed module  
 is always parity-twisted ($\epsilon=-1$, Lemma \ref{lemfixed} (i)), whereas in 
 Example \ref{exB} ($n=2$) the fixed twisted module $T_1$ has $\epsilon=1$ and $T_2$ has
  $\epsilon=-1$. 
\end{remark}

\section{The twisted Verlinde formula for $V$}
\label{secmain}

Let $\Q^{\rm fix}=\{Q\in\Q:J\fus Q\simeq Q\}$. 
It is empty if $d=-1$ by Lemma \ref{lemfree}).  
 Let $\Q^{\rm free}/J$ be the set of two-element orbits $[Q]=\{Q,J\fus Q\}$ of $J$ on $\Q$. 
 By Lemma \ref{lemvals} (4) $\eta$ is constant on orbits. 
 Set $\Q^{\rm free}_\pm=\{[Q]:\eta_Q=\pm1\}$ and  
$\Q^{\rm fix}_\pm=\{Q\in\Q^{\rm fix}:\eta_Q=\pm1\}$. By Lemma \ref{lemfixed} (i),  
$\Q^{\rm fix}_+=\left\{\emptyset\right\}$. The next main result of this paper, the 
twisted Verlinde formula for $V$, is the content of 

\begin{theorem} 
\label{thmmain}
 Adopt assumptions of Subsection \ref{secsetting}. 
Let $g_1$, $g_2\in G$ with $\langle g_1,g_2\rangle=G$, $g_3=g_1g_2$, and let 
$X_i\in\M(g_i)$ be simple objects of $\Rep(\Vev)$ (free or $J$-fixed). Set 
\[
N^\pm:=N^{X_3}_{X_1,X_2}\pm N^{J\fus X_3}_{X_1,X_2}. 
\]
Let $f_0(Q)=S_{X_1,Q}S_{X_2,Q}\ol{S_{X_3,Q}}/S_{V_\ev,Q}$. Then 

\begin{eqnarray}
\label{eqmain}
(a)  \qquad \qquad \qquad \qquad 
N^\pm=\delta_{\,\pi(X_1)+\pi(X_2),\,\pi(X_3)}
\Bigl(4\sum_{[Q]\in\Q^{\rm free}_\pm}f_0(Q)\;+\;2\sum_{Q\in\Q^{\rm fix}_\pm}f_0(Q)\Bigr),
\qquad \qquad &&
\end{eqnarray}
where the Kronecker delta is taken in $\Z/2\Z$ 
(equivalently, it is one iff $\epsilon_{X_1}\epsilon_{X_2}\epsilon_{X_3}=1$), 
and $Q$ is any representative of $[Q]$.

(b) in conditions \ref{hypdict}  
$N^+$ is the dimension and $N^-$ the superdimension of the graded space 
$I(F_V(X_3);F_V(X_1),F_V(X_2))$ of intertwining operators 
(gradings by the even representatives $X_i$). $F_V(X_i)$ is $g_i$-twisted for 
$\pi(X_i)=0$ and $g_iP$-twisted for $\pi(X_i)=1$. 

\end{theorem}

Note that all fusion coefficients and the tensor product $\fus$ here are those of 
$\Rep(\Vev)$, Definition \ref{deffusion}.
  $N^{X_3}_{X_1,X_2}$ is the individual coefficient if at least two of 
$X_1$, $X_2$, $X_3$ are $G$-stable, and otherwise the orbit average 
$|O_{X_1}|^{-1}|O_{X_2}|^{-1}\sum_{N_i\in O_{X_i}}N^{X_3}_{N_1,N_2}$ of 
Theorem \ref{thmTV} (i). The same convention applies to $N^{J\fus X_3}_{X_1,X_2}$, 
and $J\fus X_3$ has the same stabilizer as $X_3$ by Lemma \ref{lemshift}. 
  The set $\Q$ refers to $G=\langle g_1,g_2\rangle$, while the entries $S_{X,Q}$ 
depend only on $g_X$ and on the normalization of $Q$, Lemma \ref{lemorbinv} (b). 

When in \eqref{eqmain} the delta equals one the summand $f_0$ 
is constant on each free orbit. Moreover, $N^{X_3}_{X_1,X_2}=\frac12(N^++N^-)$, 
$N^{J\fus X_3}_{X_1,X_2}=\frac12(N^+-N^-)$. 
 Both vanish if $\pi(X_3)\ne\pi(X_1)+\pi(X_2)$ in $\Z/2\Z$. 
 For $d=-1$ no fixed points occur and \eqref{eqmain} 
 is $4\sum_{[Q]\in\Q^{\rm free}_\pm}f_0(Q)$. 
 For $d=1$ the fixed points contribute only to $N^-$.

\begin{proof}
\emph{(a)\; 1.} By Theorem \ref{thmTV} applied to $(X_1,X_2,X_3)$ and to 
$(X_1,X_2,J\fus X_3)$ (the same stability pattern), and by Lemma \ref{lemrow} 
($\eta_Q$ is real), 
\[
N^{X_3}_{X_1,X_2}=\sum_{Q\in\Q}f_0(Q),\qquad N^{J\fus X_3}_{X_1,X_2}=\sum_{Q\in\Q}\eta_Qf_0(Q),
\]
thus $N^\pm=\sum_Q(1\pm\eta_Q)f_0(Q)=2\sum_{Q:\eta_Q=\pm1}f_0(Q)$, 
since $1\pm\eta_Q\in\{0,2\}$. This is the origin of the outer factor $2$. 

\emph{2.} Let $Q\in\Q$ and give $J\fus Q$ the induced normalization since 
$f_0$ does not depend on normalizations, Theorem \ref{thmTV}. 
  By Lemma \ref{lemcol} and $\kappa=d\epsilon$, $d^2=1$, and 
\[
f_0(J\fus Q)=\frac{\kappa_{X_1}\kappa_{X_2}\kappa_{X_3}}{d}\,f_0(Q)=\epsilon_{X_1}\epsilon_{X_2}\epsilon_{X_3}\,f_0(Q).
\]
 In an arbitrary normalization of $J\fus Q$ the extra factor is 
 $\lambda(g_1)\lambda(g_2)\ol{\lambda(g_3)}=\lambda(g_1g_2g_3^{-1})=1$.  If $Q$ is fixed, 
 $J\fus Q\simeq Q$ and $f_0(J\fus Q)=f_0(Q)$, thus $f_0(Q)=0$ unless 
 $\epsilon_{X_1}\epsilon_{X_2}\epsilon_{X_3}=1$. 

\emph{3. Orbit decomposition:} 
  $\Q$ is the disjoint union of $\Q^{\rm fix}$ and the free orbits 
  $[Q]=\{Q,J\fus Q\}$ of size two ($\Q^{\rm fix}=\left\{\emptyset\right\}$ 
  if $d=-1$, so then all orbits are free). For a free orbit, 
  $\eta_{J\fus Q}=\eta_Q$ (see Lemma \ref{lemvals} (4)). 
 Thus, both members belong to the same set $\{\eta=1\}$ or $\{\eta=-1\}$, 
  and, according to part 2. above, their contributions are $f_0(Q)$ and 
  $\epsilon_{X_1}\epsilon_{X_2}\epsilon_{X_3}f_0(Q)$. 
  Let $\delta=1$ if $\epsilon_{X_1}\epsilon_{X_2}\epsilon_{X_3}=1$ and $\delta=0$ otherwise. 
  Then
\[
\sum_{\eta_Q=\pm1}f_0(Q)=\sum_{[Q]\in\Q^{\rm free}_\pm}\bigl(f_0(Q)+f_0(J\fus Q)\bigr)+\sum_{Q\in\Q^{\rm fix}_\pm}f_0(Q)=\delta\Bigl(2\sum_{\Q^{\rm free}_\pm}f_0+\sum_{\Q^{\rm fix}_\pm}f_0\Bigr),
\]
since 
$f_0(Q)+f_0(J\fus Q)=(1+\epsilon_{X_1}\epsilon_{X_2}\epsilon_{X_3})f_0(Q)=2\delta f_0(Q)$, 
 and $f_0(Q)=0$ at a fixed point when $\delta=0$ (see part 2). 
 The inner $2$ is the size of a free $J$-orbit (two equal terms), 
 a fixed point is a single term. Multiplying by the outer $2$ of part 1. above 
 gives the coefficients $4$ (free orbits) and $2$ (fixed points) in \eqref{eqmain}.
  Since $\epsilon_{X_1}\epsilon_{X_2}\epsilon_{X_3}=1$ iff $\pi(X_1)+\pi(X_2)=\pi(X_3)$ 
 in $\Z/2\Z$, this is the delta of \eqref{eqmain}. 

(b): By the second condition of \ref{hypdict}
 the even and odd parts of $I(F_V(X_3);F_V(X_1),F_V(X_2))$ 
have dimensions $N^{X_3}_{X_1,X_2}$ and $N^{J\fus X_3}_{X_1,X_2}$. 
 Their sum and difference are $N^\pm$. The twist statement 
 is given by the first condition of \ref{hypdict} and Lemma \ref{lemGP}. 
\end{proof}

We next adopt the following parity conventions and dependence on representatives. 
Replacing the even representative $X_1$ by $J\fus X_1$ (i.e., passing from $F_V(X_1)$ to 
$\Pi F_V(X_1)$) gives, using $J\fus J\simeq V_\ev$, 
$N^{X_3}_{J\fus X_1,X_2}=N^{J\fus X_3}_{X_1,X_2}$ and 
$N^{J\fus X_3}_{J\fus X_1,X_2}=N^{X_3}_{X_1,X_2}$. 
 Hence, $N^+\mapsto N^+$ and $N^-\mapsto-N^-$. This agrees with \eqref{eqmain}, 
since $f_0$ changes by the factor $\eta_Q$ (Lemma \ref{lemrow}), which is $1$ on the 
domain of $N^+$ and $-1$ on that of $N^-$. Thus, the dimension is independent of the 
parity conventions and the superdimension changes sign with the parity of an odd number 
of inputs/outputs. The next of this paper main results concerns the fusion product. Let us 
work in conditions \ref{hypdict} and \ref{hypmono}. 

\begin{corollary}
\label{corproduct}

Let $X_1$, $X_2$ be $J$-free with $\langle g_1,g_2\rangle=G$, and assume that every 
simple constituent of $X_1\fus X_2$ is $J$-free (if $g_3$ generates $G$ this can be tested 
on the $S$-matrix by Lemma \ref{lemfixed} (iii)). Then in $\Rep(V)^{\rm free}$
\[
F_V(X_1)\fus_VF_V(X_2)\simeq\bigoplus_{[X_3]}
\Bigl(N^{X_3}_{X_1,X_2}F_V(X_3)\oplus N^{J\fus X_3}_{X_1,X_2}\,\Pi F_V(X_3)\Bigr), 
\]
the sum over the $J$-orbits $[X_3]=\{X_3,J\fus X_3\}$ of simple constituents with 
$g_3=g_1g_2$, $\pi(X_3)=\pi(X_1)+\pi(X_2)$ (in $\Z/2\Z$), and multiplicities are given 
by \eqref{eqmain}.
\end{corollary}

\begin{proof}
  The second conditions of \ref{hypmono} gives $F_V(X_1)\fus_VF_V(X_2)\simeq F_V(X_1\fus X_2)
=\bigoplus_{X_3}N^{X_3}_{X_1,X_2}F_V(X_3)$, the sum over all simple $X_3$ 
(both members of each $J$-orbit). Let us group the terms by orbits. Since  
$F_V(J\fus X_3)=\Pi F_V(X_3)$ and $N^{J\fus X_3}_{X_1,X_2}$ is the multiplicity of 
the second member, the orbit $[X_3]$ contributes 
$N^{X_3}_{X_1,X_2}F_V(X_3)\oplus N^{J\fus X_3}_{X_1,X_2}\Pi F_V(X_3)$. 
  Replacing the chosen representative $X_3$ by $J\fus X_3$ interchanges the two summands. 
 I.e., using $J\fus J\simeq\Vev$, $N^{J\fus X_3}_{X_1,X_2}$ and $N^{X_3}_{X_1,X_2}$ swap and 
$F_V(J\fus X_3)=\Pi F_V(X_3)$, $\Pi\Pi=\id$. Thus, the displayed sum does not depend on 
the representatives. Vanishing for $\pi(X_3)\ne\pi(X_1)+\pi(X_2)$ is given by Theorem 
\ref{thmmain} (a).
\end{proof}

\section{Consistency checks}
\label{secchecks}

In this Section we do some consistency checks. 

\begin{proposition}
\label{propvac}
Let $g_2$ generate $G$ and $X_2$, $X_3\in\M(g_2)$. For $X_1=V_\ev$: 
$N^+=\delta_{X_3,X_2}+\delta_{X_3,J\fus X_2}$ and 
$N^-=\delta_{X_3,X_2}-\delta_{X_3,J\fus X_2}$, as required by 
$V\fus_VF_V(X_2)\simeq F_V(X_2)$.
\end{proposition}
\begin{proof}
Here $f_0(Q)=S_{X_2,Q}\ol{S_{X_3,Q}}$ and part 1 of the proof of Theorem \ref{thmmain} 
gives $N^\pm=\sum_Q(1\pm\eta_Q)S_{X_2,Q}\ol{S_{X_3,Q}}$. 
By Proposition \ref{propunit}, $\sum_QS_{X_2,Q}\ol{S_{X_3,Q}}=\delta_{X_2,X_3}$, 
and by Lemma \ref{lemrow} and Proposition \ref{propunit} (note that $J\fus X_2\in\M(g_2)$),  
$\sum_Q\eta_QS_{X_2,Q}\ol{S_{X_3,Q}}
=\sum_QS_{J\fus X_2,Q}\ol{S_{X_3,Q}}=\delta_{J\fus X_2,X_3}$.
\end{proof}
The selection rules are decribed by 

\begin{corollary}
\label{corselect}
Under the assumptions of Theorem \ref{thmmain} (a): 

(1) $N^{X_3}_{X_1,X_2}=N^{J\fus X_3}_{X_1,X_2}=0$ if $\pi(X_3)\ne \pi(X_1)+\pi(X_2)$ 
in $\Z/2\Z$; 

(2) otherwise $N^+\ge|N^-|$ and $N^+\equiv N^-\pmod2$, for individual coefficients.
\end{corollary}

The dimension sum rule is given by 

\begin{proposition} 
\label{propsum}
Let $g_1$, $g_2\in G$ and assume that $g_3=g_1g_2$ generates $G$ 
(then $\langle g_1,g_2\rangle=G$, thus Theorem \ref{thmTV} applies). 
Set $\delta(X)=S_{X,V_\ev}/S_{V_\ev,V_\ev}$. For $X_1\in\M(g_1)$, $X_2\in\M(g_2)$, 
\[
\sum_{X_3\in\M(g_3)}N^{X_3}_{X_1,X_2}\,\delta(X_3)=\delta(X_1)\,\delta(X_2),
\]
where $N$ is replaced by the orbit average of Theorem \ref{thmTV} (i) if neither 
$X_1$ nor $X_2$ is stable.
\end{proposition}

\begin{proof}
All $X_3\in\M(g_3)$ are stable, thus by Theorem \ref{thmTV} the left side equals to 
$\sum_{X_3}\mathsf T(X_1,X_2,X_3)\delta(X_3)
=\sum_Q\frac{S_{X_1,Q}S_{X_2,Q}}{S_{V_\ev,Q}}\cdot\frac{1}{S_{V_\ev,V_\ev}}
\sum_{X_3}\ol{S_{X_3,Q}}S_{X_3,V_\ev}$. By Proposition \ref{propunit} (columns, 
for the generating sector $g_3$) the inner sum is $\delta_{Q,V_\ev}$, which leaves 
$S_{X_1,V_\ev}S_{X_2,V_\ev}/S_{V_\ev,V_\ev}^2$.
\end{proof}

Since $\delta(X)\cdot|O_X|$ is the dimension of the $U$-module $X_\chi$ 
(Lemma \ref{lemorbrel} (b)), this is the usual dimension count of the fusion rules. 
 It is an independent consistency check of statements from Lemma \ref{lemorbrel} to Theorem 
\ref{thmTV}, and it is verified in all examples (see Section \ref{secexamples}).

\section{Inner automorphisms}
\label{secinner}

 Let $h\in\Vev$ be of weight one with $L(m)h=\delta_{m,0}h$ and 
$h_mh=\delta_{m,1}\langle h,h\rangle\mathbf 1$ ($m\ge 0$).  
such that $h_0$ acts semisimply on $V$ with spectrum in $\frac1n\Z$, 
 Let $\sigma=e^{2\pi ih_0}$ have order $n$ and act faithfully on $\Vev$ as in \cite{DL}.  
 Li operator $\Delta(h,z)=z^{h_0}\exp\bigl(\sum_{m\ge1}\frac{h_m}{-m}(-z)^{-m}\bigr)$ turns 
a module $M$ into $M^{(h)}=(M,Y_M(\Delta(h,z)\cdot,z))$ \cite{Li}. 
     We consider such $h_0$ that has the following property: it acts semisimply on every 
irreducible untwisted $\Vev$-module
 (we call this property in what follows the property of Section \ref{secinner}).  

Then we have the following statement about inner automorphisms. 

\begin{proposition}
\label{propinner}
\begin{enumerate}

\item[(i)] Every irreducible $\sigma^{\pm a}$-twisted $\Vev$-module is $M^{(ah)}$ for 
 an irreducible $\Vev$-module $M$ \cite{Li, DLM1}. Moreover, when the first argument 
is untwisted, $N^{M_3^{(ah)}}_{M_1,M_2^{(ah)}}=N^{M_3}_{M_1,M_2}$ for untwisted 
$M_1$, $M_2$, $M_3$ \cite{DLM1}.

\item[(ii)] Every irreducible twisted $\Vev$-module is $G$-stable. Hence, $\Q=\M(1)$, 
and Theorem \ref{thmTV} (ii) holds for all $g_1$, $g_2$ with 
 $\langle g_1,g_2\rangle=G$ and all 
$n$ (including composite and even $n$), 
and Theorem \ref{thmmain} (a) holds for all simple $X_i$.

\item[(iii)] Adopt assumptions \ref{hypL}. Then conditions \ref{hypdict}  
and \ref{hypmono} follow from the untwisted statements of \cite{CKM, HKL}  
by transport of structure. 
\end{enumerate}
\end{proposition}

\begin{proof}
(i) is given by the corresponding citation results mentioned above with proofs using 
the Jacobi identity and $\Delta(h,z)$.  

(ii) By (i) the module is $M^{(ah)}$. The zero mode of $h$ on $M^{(ah)}$ is  
$h_0^M+a\langle h,h\rangle$ \cite{Li}, semisimple by the condition of this Section). 
 It commutes with $L(0)$, 
and satisfies $[h_0^{M^{(ah)}},Y(v,z)]=Y(h_0v,z)$ since $h$ is fixed by $\sigma^a$ and has 
integer modes. Hence, $e^{2\pi ih_0^{M^{(ah)}}}:M^{(ah)}\to M^{(ah)}\circ\sigma$ is an 
isomorphism (cf. \cite{DL}). 

(iii) is a transport of structure through $\Delta$, 
with the following dependence. The first condition of \ref{hypdict}: 
 by the first condition of \ref{hypL}  $F_V(M^{(ah)})=F_V(M)^{(ah)}$ is a module 
twisted by $\phi(ah)$ with the same parity twist as $F_V(M)$, and $F_V(M)$ is an 
(untwisted or parity-twisted) $V$-module by the untwisted theory \cite{CKM}. 
 The second condition of \ref{hypdict}: for $X_3$ twisted and  one of $X_1$, $X_2$ untwisted, 
 the second conditions of \ref{hypL} transports the graded space of intertwining operators of the untwisted case. 
 For $X_1$, $X_2$ both twisted one requires the third condition of \ref{hypL})     
 (including its compatibility with the $\Z/2\Z$-gradings, 
 which is part of the formulation of the third condition of \ref{hypL}).  
  Conditions \ref{hypmono}: the classification and the tensor product on $\Rep(V)^{\rm free}$ for 
twisted sectors are transported from the untwisted ones by conditions \ref{hypL} 
 if one additionally 
assumes that the tensor product of twisted $V$-modules exists and that $M\mapsto M^{(ah)}$ 
is compatible with it. 
\end{proof}

\begin{remark}
\label{remtwotw}
Let us discuss (i).  
  The results \cite{Li, DLM1} twist the second and third arguments of an intertwining 
operator and leave the first untwisted. A fusion coefficient of two twisted modules,  
$N^{M_3^{((a+b)h)}}_{M_1^{(ah)},M_2^{(bh)}}$, is not of this form. The statement that it 
equals $N^{M_3}_{M_1,M_2}$ is the third condition of \ref{hypL}, which we do not prove. A proof would use that 
$\Vev^{(ah)}$ is an invertible object for a tensor product of twisted modules and the 
associativity of intertwining operators, i.e., a piece of the $G$-crossed theory.  
  The numerical evidence below (all examples with inner automorphisms, in particular, 
the agreement of the formulas with the closed-form rule in \texttt{verify\_grouplaw.p y}, 
and Example \ref{exF}) is consistent with the third condition of \ref{hypL}. 
 By (ii) all modules are stable, 
 thus Theorem \ref{thmTV} (ii) is available for every generating triple and all coefficients 
are computed from $U$ by part 0. 
\end{remark}

In Examples \ref{exA}-\ref{exE} 
the automorphism is of this type and the condition of this Section  
 holds ($h_0$ is a lattice or Cartan zero mode).  

\section{Examples}
\label{secexamples}

\subsection{Computations in examples}

In computations related to our examples below we used Plain Python 3. 
 The code is contained in \texttt{python3 verify\_all.py} etc. 
 Outputs are is given in \texttt{verification\_log.txt}. 

Examples \ref{exA}-\ref{exE} is in \texttt{verify\_all.py}. 
 Independent test of Lemma \ref{lemorbrel} in Example \ref{exA} from twisted traces 
(see Remark \ref{remconv}) is in \texttt{verify\_twisted\_traces.py}. 

$(ST)^3=S^2$ and the Gauss identity of Lemma \ref{lemgauss} 
on the modular data used here ($V_{\Z\beta}$ with $\beta^2=4,6,18$; Virasoro $M(3,5)$; 
Ising$\,\boxtimes\,\Z/4$); residuals $\le4\cdot10^{-15}$ are in \texttt{verify\_gauss.py}. 

Comparison of the formulas (Theorem \ref{thmTV}, 
\eqref{eqmain}) with fusion coefficients given by the closed form spectral flow rule 
in Examples \ref{exA}, \ref{exC}, \ref{exD} (see Remark \ref{remtwotw}) is in 
\texttt{verify\_grouplaw.py}. 

Example \ref{exF} is in \texttt{verify\_product.py}. 

The barred/unbarred illustration in Remark \ref{remconv} is in \texttt{verify\_conv.py}. 

Each script header lists its inputs and its normalization conventions. 
 These are the conventions of Section \ref{secS} 
 (in particular $\varphi_{Q_k}=e^{-2\pi ik/2n}\sigma$ in Example \ref{exA}, 
 $e^{-2\pi ik/n}\sigma$ in Example \ref{exC}). 
 The scripts are supplementary material and are not reproduced here.
 
\emph{Inputs and checks:} the inputs to \texttt{verify\_all.py} are 
 
(i) the modular data of $U$ ($S^U$ and the fusion rules of $\mathcal C_U$); 

(ii) for every $\Vev$-module the list of its $U$-components 
(Proposition \ref{factclass}, each list is justified in the example by counting). 
Other subjects are computed: 

\noindent(R) Lemma \ref{lemorbrel} (b) is checked for \emph{all} components of all modules. 
 The orbifold $S$-matrix is computed as $S_{X,Q}=|G_X|S^U_{X_{\chi_0},Q_{t=0}}$.

\noindent(U) Proposition \ref{propunit}: the block $(S_{X,Q})_{X\in\M(g),Q\in\Q}$ 
 is square and unitary for every generating sector $g$.

\noindent(F) Fusion coefficients are computed by Part 0 of Theorem \ref{thmTV} from the fusion 
rules of $U$ ($U$-components, divided by $|O_1||O_2|$). 

\noindent(TV, 4.1) Theorem \ref{thmTV} and formula \eqref{eqmain} 
(with the fixed-point weights)
 are verified on all in-scope triples ($\langle g_1,g_2\rangle=G$, $g_3=g_1g_2$). 
 Where all modules are stable, the triples with $\langle g_1,g_2\rangle\ne G$ 
 (including $g_1=g_2=1$), which lie outside Theorems \ref{thmTV} and \ref{thmmain} as stated, 
 were checked too and are counted separately (``extra'').

\noindent(SR, W, 3.9) the sum rule (Proposition \ref{propsum}); the weight criterion of Lemma 
\ref{lemweight} and Remark \ref{remweights}; Lemma \ref{lemfixed} including the criterion 
(iii).

\begin{center}\small
\begin{tabular}{llrrrrrrr}
\hline
Example & $d$ & (R) & unitary blocks & in-scope & extra & SR & W & 3.9(iii)\\ \hline
\ref{exA}, $n=3$ & $1$ & 432 & 2$\times$($4\times4$) & 512 & 64 & 96 & 12 & 8\\
\ref{exA}, $n=9$ & $1$ & 11664 & 6$\times$($4\times4$) & 4608 & 576 & 864 & 36 & 24\\
\ref{exA}, $n=15$ & $1$ & 54000 & 8$\times$($4\times4$) & 12288 & 2112 & 1920 & 60 & 32\\
\ref{exB}, $p=2$ & $1$ & 40 & 1$\times$($2\times2$) & 48 & - & 16 & 4$^*$ & 2\\
\ref{exC}, $n=3$ & $-1$ & 1728 & 2$\times$($8\times8$) & 4096 & 512 & 384 & 24 & 16\\
\ref{exC}, $n=9$ & $-1$ & 46656 & 6$\times$($8\times8$) & 36864 & 4608 & 3456 & 72 & 48\\
\ref{exC}, $n=15$ & $-1$ & 216000 & 8$\times$($8\times8$) & 98304 & 16896 & 7680 & 120 & 64\\
\ref{exD}, $n=4$ & $1$ & 1024 & 2$\times$($4\times4$) & 768 & 256 & 128 & 4$^*$ & 8\\
\ref{exE}, $p=2$ & $1$ & 72 & 1$\times$($3\times3$) & 81 & 27 & 18 & 3$^*$ & 3\\
\ref{exF}, $n=4$ & $1$ & - & 2$\times$($4\times4$) & 1536 & - & 256 & - & -\\
\end{tabular}

\smallskip
\footnotesize $^*$ untwisted modules only (for even $n$ the twisted weight test is uninformative). ``Extra'' triples are not covered by the theorems as stated; in Example \ref{exB} there are unstable modules and no extra triples are checked. Example \ref{exF} is run by a separate script from products of the data of its factors, thus (R), W and 3.9(iii) do not apply to it.
\end{center}

\subsection{The content of tests for each example} 

Stable-column formula (Lemma \ref{lemorbrel} (b)) was tested independently of $S^U$ 
in Example \ref{exA} with twisted traces, $n=3$, $9$.
 In Examples \ref{exA} and \ref{exC} composite $n$ are $9$ and $15$. 
In Examples \ref{exD} and \ref{exF} composite $n$ is $4$, and no composite $n$ in 
Examples \ref{exB} and \ref{exE}.  
 Non-trivial stabilizers are testes in Example \ref{exB} with $n=2$,
 and only by factorization in Example \ref{exF}. 
In every example the checks (R), (U), (F), (TV,\,4.1), (SR) test the algebra 
of the proofs and the stated conventions given the inputs. 
 They do not test Propositions \ref{factrep} and \ref{factclass}. 

\begin{example}
\label{exA}

In this example we consider complex fermion with inner automorphism of odd order $n$. 
Let $\alpha$ have $(\alpha,\alpha)=1$, $V=V_{\Z\alpha}$, $\Vev=V_L$ with $L=2\Z\alpha$, 
$J=V_{\alpha+L}$ (weights in $\Z+\frac12$, thus $d=1$), $\sigma=e^{2\pi i\alpha(0)/n}$, 
$n\in\{3,9,15\}$ (prime, respectively, composite).

\emph{Untwisted data:} $Q_k=V_{L+k\alpha/2}$, $k\in\Z/4$, $S_{Q_k,Q_l}=\frac12e^{-2\pi ikl/4}$,
 $J=Q_2$, $\eta_{Q_k}=(-1)^k$. $J$ freely acts by $k\mapsto k+2$. 
 All $Q_k$ are stable: $\Q=\{Q_0,\dots,Q_3\}$, $\Q^{\rm free}_+=\{[Q_0]\}$, 
 $\Q^{\rm free}_-=\{[Q_1]\}$. 

\emph{The fixed-point algebra:} $U=V_{L'}$, $L'=2n\Z\alpha$ (norm $4n^2$), 
 $\mathcal C_U$ is pointed, $\Z/4n^2$, modules $V_{L'+j\alpha/(2n)}$, 
 $S^U_{jj'}=(2n)^{-1}e^{-2\pi ijj'/4n^2}$. 
 The module $Q_k$ consists of the $U$-modules $j=nk+4ns$, 
 the $\sigma^a$-twisted module $T^{(a)}_m=V_{L+(m/2+a/n)\alpha}$ ($a\in\Z/n$, $m\in\Z/4$)
 of $j=nm+2a+4ns$, $s\in\Z/n$, with $T^{(a+n)}_m=T^{(a)}_{m+2}$. 
 The count $4n^2=n\cdot(n\cdot4)$ confirms Proposition \ref{factclass}. 
 $V^G=V_{n\Z\alpha}$ is a simple lattice VOSA.

\emph{Orbifold data:} with $\varphi_{Q_k}=e^{-2\pi ik/2n}\sigma$ 
 (a representation, since $\sigma^n=(-1)^k$ on $Q_k$. 
 The label of the component $j=nk+4ns$ is $-2s$), Lemma \ref{lemorbrel} gives 
\[
S_{T^{(a)}_m,Q_k}=\tfrac12\,e^{-2\pi imk/4}\,e^{-2\pi iak/2n},
\qquad \kappa_{T^{(a)}_m}=(-1)^m=\epsilon_{T^{(a)}_m}.
\]
 even $m$ corresponds to the $\sigma^a$-twisted modules $V_{\Z\alpha+a\alpha/n}$ (NS), 
 odd $m$ to the $\sigma^aP$-twisted ones (Ramond). 
  This labelling is confirmed independently by the weight test of Remark \ref{remweights}, 
 $e^{2\pi in(h_{J\fus X}-h_X)}=-d\epsilon_X$, valid for all modules and sectors. 
 $J$ acts by $m\mapsto m+2$ on all sectors. 

\emph{Independent test of the orbifold $S$-matrix:} 
 the entries above were derived from $S^U$ by Lemma \ref{lemorbrel}.
 For this family the twisted characters are theta functions and the law 
 \eqref{eqlaw} can be tested: for $n=3$, $9$, all sectors and three values of 
 $\tau$, $F_{a,m}(-1/\tau)=\sum_kS_{T^{(a)}_m,Q_k}\,G_{a,k}(\tau)$ 
 holds with residual $\le1.2\cdot10^{-15}$, where $F_{a,m}$ is the character 
 of $T^{(a)}_m$ and $G_{a,k}$ the $\sigma^{-a}$-trace on $Q_k$ 
 in the normalization $\varphi_{Q_k}=e^{-2\pi ik/2n}\sigma$ 
 (\texttt{verify\_twisted\_traces.py}). 
 Changing the sign of $a$ or conjugating $S$ gives residuals $0.68$-$1.22$. 
 This does not use the modular data of $U$, and checks Lemma \ref{lemorbrel} 
 and the conventions of Remark \ref{remconv} for this family 
 (it tests characters, i.e., $v=\mathbf 1$, and one family only).

\emph{Checks:} the twisted Verlinde formula reproduces the fusion rules 
 $T^{(a)}_m\fus T^{(b)}_{m'}=T^{(a+b)}_{m+m'}$ obtained from part 0: for $a+b<n$ the phases 
 $e^{-2\pi ik(a+b-c)/2n}$ cancel. For $a+b\ge n$ they produce $(-1)^k$, i.e., the shift 
 $m\mapsto m+2$. \eqref{eqmain} gives $N^+=4f_0(Q_0)=1$ and $N^-=4f_0(Q_1)=\pm 1$ if 
 $\pi(X_1)+\pi(X_2)=\pi(X_3)$ (and zero otherwise). 
   
\emph{Counts:} for $n=9$  
 the sectors $a\in\{3,6\}$ generate a subgroup $H$ of order three. There the theorem is to  
 be applied with $H$, and by Lemma \ref{lemorbinv} (b) the entries are those of the case 
 $n=3$ checked for each entry by the script. Here $\Q_H=\Q_G$ since all modules are stable, 
 thus the formulas agree.

Summarising the results we have the following.  
 With the first assumption of Subsection \ref{secsetting},
 $V_{\Z\alpha}$ is a simple lattice VOSA \cite{FLM}.
 $V_L$, $L$ even positive definite is simple, rational, $C_2$-cofinite, self-dual \cite{D93}. 
With the second assumption of Subsection \ref{secsetting},
 $\sigma$ acts on $e^{2m\alpha}$ by $e^{4\pi im/n}$: 
order $n$ on $\Vev$ ($n$ odd). 
 For $V^G$ simple, then one has Lemma \ref{lemfixedsimple}, in fact $V^G=V_{n\Z\alpha}$. 
With Theorem \ref{thmTV} (ii),  $n=3$ prime, $n=9$, $15$: all modules stable,  
 WIth Proposition \ref{propinner} (ii), the property of Section \ref{secinner}  
 holds for lattice zero modes.  
For fixed points, one has that $J$ acts freely on all of $\R$. 
Propositions \ref{propunit}, \ref{propsum} and  Lemma \ref{lemweight} are checked.   
\end{example}

\begin{example} 
\label{exB}
 In this example we consider two Majorana fermions, $p=2$: non-inner, unstable modules, and 
fixed twisted modules. 
  Let $V=F\otimes F$ ($F_\ev=I_0$, $F_\od=I_{1/2}$ Ising modules), 
$\Vev=I_0\otimes I_0\oplus I_{1/2}\otimes I_{1/2}$, 
$J=I_0\otimes I_{1/2}\oplus I_{1/2}\otimes I_0$, $\sigma=1\otimes P$, order two, nontrivial 
on $\Vev$.
  Then $U=I_0\otimes I_0$, $\mathcal C_U$ the product of two Ising categories 
(modules $(a,b)$, $a$, $b\in\{1,\sigma,\psi\}$), $J_0=(\psi,1)$, $d=1$. 
 The $U$-components are: $V_\ev=(1,1)+(\psi,\psi)$, $J=(\psi,1)+(1,\psi)$, 
 $T_1=(1,\sigma)+(\psi,\sigma)$, $T_2=(\sigma,1)+(\sigma,\psi)$, and the unstable pair 
 $Q_{[1]},Q_{[3]}$ with the single $U$-module $(\sigma,\sigma)$ (count $9=2\cdot4+1$).
  Thus, $\Q=\{V_\ev,J\}$, $\M(\sigma)=\{T_1,T_2\}$. Lemma \ref{lemorbrel} gives, 
 in the columns $(V_\ev,J)$,
\[
S_{T_1,\cdot}=(\tfrac{\sqrt2}2,\tfrac{\sqrt2}2),\quad S_{T_2,\cdot}
=(\tfrac{\sqrt2}2,-\tfrac{\sqrt2}2),\quad S_{Q_{[1]},\cdot}=(\tfrac12,-\tfrac12).
\]
Part 0 gives $T_1\fus T_1=T_2\fus T_2=V_\ev\oplus J$, $T_1\fus T_2=Q_{[1]}\oplus Q_{[3]}$ 
(an unstable output) and $Q_{[1]}\fus T_1=T_2$ (an unstable input). 
Theorem \ref{thmTV} reproduces all of them (48 in-scope triples, including these). 

\emph{Signs:} $\kappa_{T_1}=1$, $\kappa_{T_2}=-1$, $\kappa_{Q_{[1]}}=-1$.
 $T_1$ is the $\sigma$-twisted NS$\otimes$R module and $T_2$ the 
 $\sigma P$-twisted R$\otimes$NS module. By the physical identification, 
 the first conditions of \ref{hypdict} is the statement that this labelling is correct. 
 The weight test holds for the four untwisted modules. 
 \eqref{eqmain} holds, e.g., $X_1=T_1,X_2=T_2,X_3=Q_{[1]}$ gives $N^+=4f_0(V_\ev)=2$,  
 $N^-=0$. 

$J\fus T_i\simeq T_i$: every twisted module is a fixed point 
(Lemma \ref{lemfixed} (iii) detects this on the $S$-matrix), 
thus no triple with $g_i=\sigma$ consists of free objects and Corollary \ref{corproduct} 
does not apply. 
 Conditions \ref{hypdict} are claimed for fixed objects (conditions \ref{hypdict})  
but nothing is known about the structure of $F_V(T_i)=T_i\oplus T_i$. 
Theorem \ref{thmmain} (a) is unaffected. Note that $T_1$ is parity-untwisted 
although fixed (see Remark \ref{remweights}). 

\emph{Higher rank (not verified):} for $V=F^{\otimes2N}$ with $\sigma$ flipping one 
fermion one has $\Vev=\mathfrak{so}(2N)_1$ and $U=F^{\otimes(2N-1)}_\ev\otimes F_\ev$. 
By the standard description of $\mathfrak{so}(m)_1$ ($m$ odd) its modular category has 
the Ising $S$-matrix. Thus, the $S$-data of this example would be the same for every $N$.
 We verified only $N=1$ by script. The identification for $N\ge2$ is the standard. 

Summarizing result, we have the following. 
  With both assumptions of Subsection \ref{secsetting}, 
 $V=F^{\otimes2}$ is simple, $\Vev$ is regular, $\sigma$ is nontrivial on $\Vev$, 
and $V^G=F\otimes I_0$.  
 With $G\cap\langle P\rangle=1$,  $\sigma\ne P=P\otimes P$. 
With Theorem \ref{thmTV}, $n=2$ prime \cite{DL}, 
 unstable $Q_{[1]}$ as input and output (part 0 is checked). 
 Propositions \ref{propunit}, \ref{propsum} are checked.  
\end{example}

\begin{example}
\label{exC}
In this example we consider $L_1(\mathfrak{osp}(1|2))$:  wrong statistics, non-pointed, 
orders $3$, $9$, $15$. 

\emph{Known from the literature:}

 (a) $V=L_1(\mathfrak{osp}(1|2))$ is a simple VOSA with $\Vev=V_{\Z\alpha}\otimes M(3,5)$,
 $(\alpha,\alpha)=2$ ($V_{\Z\alpha}=L_1(\mathfrak{sl}_2)$) 
 and $J=V_{\Z\alpha+\alpha/2}\otimes\phi_{1,4}$,  
  the coset construction of \cite{CFK}.   
  One requires that $V_{\Z\alpha+\alpha/2}\otimes\phi_{1,4}$ is a simple current of order 
 two for $\Vev$ whose superalgebra extension $V$ is simple.   

(b) the Virasoro minimal model $M(3,5)$ is rational and $C_2$-cofinite, \cite{W}. 

\emph{We checked here:} the central charge $\frac25=1+(-\frac35)$ and the weight 
 $\frac14+\frac34=1$ of $J$. The $S$-matrix 
 $S_{\mu\nu}=\sqrt{2/5}\,(-1)^{\mu+\nu}\sin(3\pi\mu\nu/5)$ of $M(3,5)$ is the modular 
 transformation of the Rocha-Caridi characters $\chi_{1,\mu}$ (residual $<10^{-15}$ at four 
 values of $\tau$), is symmetric and orthogonal, and gives non-negative integer fusion 
 coefficients with $\phi_{1,4}\fus\phi_{1,4}=\phi_{1,1}$, 
 $\phi_{1,2}\fus\phi_{1,2}=\phi_{1,1}\oplus\phi_{1,3}$, 
 $\dim_q\phi_{1,\mu}=(1,\,0.618,\,-0.618,\,-1)$. 
 Hence, $\dim_qJ=1\cdot(-1)$, $\theta_J=e^{2\pi i}=1$: $d=-1$ and $\theta_J=-d$ 
 (Proposition \ref{factJ} (ii)). 

Everything below in this example, in particular the value $d=-1$, the sign data and 
counts depend on input (a). 
  Take $\sigma=\sigma_n=e^{2\pi i\alpha(0)/n}\otimes1$, 
$n\in\{3,9,15\}$: order $n$ on $\Vev$ and on $J$, trivial on the Virasoro factor.

\emph{Data:} modules of $\Vev$: $(\ell,\mu)$, $\ell\in\Z/2$. 
 $J=(1,4)$ freely acts by $(\ell,\mu)\mapsto(\ell+1,5-\mu)$, Lemma \ref{lemfree}. 
 Twisted modules: $(T^{(a)}_\ell,\mu)$, $T^{(a)}_\ell=V_{\Z\alpha+(\ell/2+a/n)\alpha}$. 
 All modules are stable (inner). The lattice factor has $U_{\rm lat}=V_{n\Z\alpha}$, 
 $\mathcal C_U$ of rank $2n^2$, $S^U_{jj'}=\frac1{n\sqrt2}e^{-2\pi ijj'/2n^2}$. 
 The $U$-components of $(T^{(a)}_\ell,\mu)$ are $(j=n\ell+2a+2ns,\mu)$, $s\in\Z/n$. 
 In the normalization $\varphi=\sigma$ the orbifold entries are the sector independent numbers
\[
S_{(T^{(a)}_\ell,\mu),(Q_k,\nu)}=\tfrac1{\sqrt2}(-1)^{\ell k}\,S_{\mu\nu}. 
\]
 Here the script uses the normalization with $\varphi_{Q_k}=e^{-2\pi ik/n}\sigma$, 
which multiplies this by $e^{-2\pi iak/n}$ in the sector $a$ by \eqref{eqnormchange}. 

\emph{Sign data:} $\eta_{(k,\nu)}=(-1)^{k+\nu}$, $d=\eta_{(0,1)}=-1$, 
$\epsilon_{(\ell,\mu)}=-(-1)^{\ell+\mu}$. The NS ($\epsilon=1$) $J$-orbits are 
$\{(0,1),(1,4)\}$, $\{(0,3),(1,2)\}$; the Ramond ones $\{(0,2),(1,3)\}$, $\{(0,4),(1,1)\}$. 
As $d=-1$, $\eta=1$  on the Ramond orbits: $N^+$ is a sum over parity-twisted and $N^-$ 
over parity-untwisted stable modules, as in \cite{CKM}. The weight test of Remark 
\ref{remweights} holds for all modules and sectors.

\emph{Checks:} the fusion coefficients computed by Part 0 are 
$N^{(c,\ell_3,\mu_3)}_{(a,\ell_1,\mu_1),(b,\ell_2,\mu_2)}
=\delta_{c,a+b}\delta_{\ell_3,\ell_1+\ell_2}N^{\mu_3}_{\mu_1\mu_2}$ 
with the (non-pointed) Virasoro fusion coefficients; Theorem \ref{thmTV} and \eqref{eqmain} 
reproduce them (counts in the table). Example: $X_1=X_2=(T^{(1)}_0,\phi_{1,2})$ (both Ramond) 
gives $X_3=(T^{(2)}_0,\phi_{1,1})$ or $(T^{(2)}_0,\phi_{1,3})$ (NS), $N^\pm=1$, consistent 
with $\epsilon_{X_1}\epsilon_{X_2}\epsilon_{X_3}=1$.

Summarizing results we have the following. 
 With the first assumption of subsection \ref{secsetting}, inputs (a), (b) above and    
 consistency are checked.   
 With the second assumption of Subsection \ref{secsetting},  
 $\sigma$ is inner on the lattice factor, 
order $n$ on $\Vev$ ($n$ odd). $V^G$ is simple by Lemma \ref{lemfixedsimple}. 
With Theorem \ref{thmTV} (ii),  $n=3$ prime; for $n=9$, $15$: all modules are stable.  
 Fixed points: none ($d=-1$, Lemma \ref{lemfree}. 
 Propositions \ref{propunit}, \ref{propsum} are checked with $8\times8$ blocks.  
This is the only example with $d=-1$ and it is based on (a). 
It shows that \eqref{eqmain} carries information beyond the twisted fusion coefficients, 
namely the sign restriction of the sums.
\end{example}

\begin{example}
\label{exD}
In this example of rank two, even composite order $n=4$. 
 Let $V=V_{\Z^2}$ (two complex fermions, $(\alpha_i,\alpha_j)=\delta_{ij}$), 
$\sigma=e^{2\pi i\alpha_1(0)/4}$. Then $\Vev=V_{L_0}$, 
 $L_0=\{(k_1,k_2):k_1+k_2\ \text{even}\}$, $J=V_{L_0+(1,0)}$, $d=1$; 
$\sigma$ has order $4$ on 
$V$ and on $\Vev$, with $\sigma^2=(-1)^{k_1}\neq P=(-1)^{k_1+k_2}$, 
thus $G\cap\langle P\rangle=1$ and $G\times\langle P\rangle\simeq\Z/4\times\Z/2$ 
is not cyclic. 
  The category of $\Vev$ is pointed, $\Z/2\times\Z/2$: $V_\ev$, $J$, 
$A=(\frac12,\frac12)$, $B=(\frac12,-\frac12)$ ($J\fus A=B$); $\eta=(+,+,-,-)$. 
$U=V_{4\Z\oplus2\Z}$ (rank $64=n^2\cdot4$), modules $(j_1,j_2)\in\Z/16\times\Z/4$ with 
$S^U=\frac18e^{-2\pi i(j_1j_1'/16+j_2j_2'/4)}$; the $U$-components of 
$T^{(a)}_Q=V_{L_0+\lambda_Q+(a/4,0)}$ are $(4\lambda_{Q,1}+a+4t,\,2\lambda_{Q,2}+2t)$, 
$t\in\Z/4$. All modules are stable (inner); $\Q^{\rm free}_+=\{[V_\ev]\}$, 
$\Q^{\rm free}_-=\{[A]\}$. 
  Theorem \ref{thmTV} reproduces the Part 0 fusion rules and 
\eqref{eqmain} holds on all in-scope triples, $\epsilon$ takes both values.
  This example has even $n$ (thus $\sigma^{n/2}\ne P$ is essential) and a non-cyclic group 
$G\times\langle P\rangle$. Two assumptions of Subsection \ref{secsetting}, 
unitarity, sum rule and the untwisted weight test are as in Example \ref{exA}/the table. 
The correspondences depend on conditions \ref{hypL}.
\end{example}

\begin{example}
\label{exE}
In this example we consider three Majorana fermions: a $J$-fixed untwisted module. 
 Let $V=F^{\otimes3}=V_{\Z\alpha}\otimes F$ ($(\alpha,\alpha)=1$), $\Vev=L_2(\mathfrak{sl}_2)$ 
($c=3/2$), $J$ the spin-one module ($h=\frac12$, $d=1$), $\sigma=e^{\pi i\alpha(0)}\otimes1$ 
(a rotation by $\pi$; order two, inner on $\Vev$). 
  Untwisted modules (all stable): 
$[0]=V_\ev$, $[1]=J$, $[\frac12]$ ($h=3/16$); $J$ swaps $[0]\leftrightarrow[1]$ 
and \emph{fixes} $[\frac12]$. $U=I_0\otimes V_{2\Z\alpha}$ (Ising$\,\boxtimes\,\Z/4$, 
$12=2\cdot6$ modules); as $U$-modules $[0]=(1,0)+(\psi,2)$, $[1]=(\psi,0)+(1,2)$, 
$[\frac12]=(\sigma,1)+(\sigma,3)$ and the three $\sigma$-twisted modules are 
$T_A=(1,1)+(\psi,3)$, $T_B=(1,3)+(\psi,1)$, $T_C=(\sigma,0)+(\sigma,2)$ 
(the components of each module are related by multiplication with $(\psi,2)$, 
the $\sigma$-odd part of $\Vev$).
  Lemma \ref{lemorbrel} gives, in the columns $([0],[1],[\frac12])$,
\[
S_{T_A,\cdot}=(\tfrac12,\tfrac12,-\tfrac{i\sqrt2}2),\quad 
S_{T_B,\cdot}=(\tfrac12,\tfrac12,\tfrac{i\sqrt2}2),\quad 
S_{T_C,\cdot}=(\tfrac{\sqrt2}2,-\tfrac{\sqrt2}2,0),
\]
and $S_{[\frac12],\cdot}=(\frac{\sqrt2}2,-\frac{\sqrt2}2,0)$, $\eta=(1,1,-1)$: the fixed point 
$[\frac12]$ has $\eta=-1$ (Lemma \ref{lemfixed} (i)); by Lemma \ref{lemfixed} (iii) the fixed 
twisted module $T_C$ is detected by $\sum_Q\eta_Q|S_{T_C,Q}|^2=\frac12+\frac12-0=1$, while 
$\sum_Q\eta_Q|S_{T_A,Q}|^2=\frac14+\frac14-\frac12=0\ne1$ ($T_A$ is free). $J\fus T_A=T_B$, 
 $J\fus T_C=T_C$; $\epsilon=(+,+,-)$ on $(T_A,T_B,T_C)$, $\epsilon_{[1/2]}=-1$.

\emph{The fixed point matters:} take $(X_1,X_2,X_3)=(T_A,T_B,[0])$: $N^{[0]}_{T_A,T_B}=1$,  
$N^{[1]}_{T_A,T_B}=0$, thus $N^+=N^-=1$. By \eqref{eqmain}: $N^+=4f_0([0])=4\cdot\frac14=1$, 
and $N^-=4\sum_{\Q^{\rm free}_-}f_0+2f_0([\tfrac12])=0+2\cdot\tfrac12=1$: the minus sign comes 
 entirely from the fixed point, with weight two (weight $4$ would give $2$, weight $0$ would 
give $0$). Over the 81 in-scope triples, treating $[\frac12]$ as a free orbit gives a wrong 
$N^-$ in $24$ cases. Also $T_C\fus T_C=[0]\oplus[1]$ ($N^+=2$, $N^-=0$) and 
$T_A\fus T_C=[\frac12]$ ($N^+=2$, $N^-=0$ since $S_{T_C,[1/2]}=0$, in accordance with Lemma 
\ref{lemfixed} (ii)).

Summarizing results, the script has checked that, 
 with both assumptions of Subsection \ref{secsetting},       
 $F^{\otimes3}$ is simple, $L_2(\mathfrak{sl}_2)$ is regular, 
 $\sigma$ is nontrivial on $\Vev$, and $V^G=U\oplus J_0=F\otimes V_{2\Z\alpha}$. 
For Theorem \ref{thmTV} (ii), $n=2$ prime, all modules are stable. 
For fixed points, $[\frac12]\in\Q^{\rm fix}$, $T_C$, Lemma \ref{lemfixed} (i)-(iii) 
is checked.  Conditions \ref{hypdict} are checked for fixed objects, 
depending on conditions \ref{hypL}. 

\end{example}

\begin{example}
\label{exF}

In this example we consider a product with a proper stabilizer, $n=4$.  

Let $A=L_1(\mathfrak{sl}_2)=V_{\Z\beta}$ ($(\beta,\beta)=2$), $\sigma_A=e^{2\pi i\beta(0)/8}$ 
(order $4$ on $A$), and let $V=A\otimes F\otimes F$, thus $\Vev=A\otimes B_0$ with $B_0$ the
 even part of Example \ref{exB}, $J=A\otimes J_B$, $d=1$, and 
$\sigma=\sigma_A\otimes(1\otimes P)$, $G=\langle\sigma\rangle\simeq\Z/4$. 
  The first assumption of Subsection 
\ref{secsetting} holds; $\sigma^2=\sigma_A^2\otimes1\ne1$ 
on $A$ and $\sigma^k$ is nontrivial 
on $B_0$ for $k$ odd, thus $G$ is faithful on $\Vev$ 
the second assumption of Subsection \ref{secsetting}, and  
$\sigma^k\ne P$ for all $k$ (Lemma \ref{lemGP}).

\emph{Structure:} twisted modules of $A\otimes B_0$ are tensor products of twisted modules, 
and twisted traces and fusion coefficients factorize. 
 All $A$-modules are stable. 
 A pair $(X_A,y)$ with $y\in\{Q_{[1]},Q_{[3]}\}$ (the unstable pair of Example \ref{exB}, 
in even sectors) has $G$-orbit of size two and stabilizer $\langle\sigma^2\rangle$. 
 The twisted modules of odd sectors have $B$-part $T_1$ or $T_2$ and are $J$-fixed. 
 In-scope triples include $(g_1,g_2,g_3)=(\sigma^2,\sigma,\sigma^3)$ with $M_1$ non-stable and 
$M_2,M_3$ stable, the situation of Remark \ref{remn4}.

\emph{Check (\texttt{verify\_product.py}):} the $S$-entries are the products of those of the 
factors (the $A$-factor from $U_A=V_{4\Z\beta}$, the $B$-factor from Example \ref{exB}, with 
the $B$-sector $a\bmod2$, Lemma \ref{lemorbinv} (b)). 
 The fusion coefficients are the products of Part 0 coefficients of the factors. 
 For $G=\Z/4$ with $\Q=\Q_A\times\Q_B$ ($4$ modules),  
the generating blocks ($a=1,3$) are $4\times4$ unitary (Proposition \ref{propunit}). 
 Theorem \ref{thmTV} and \eqref{eqmain} hold on all $1536$ in-scope triples 
 (of which $768$ have  one non-stable module, of stabilizer order two, and two stable ones; 
 all involve $J$-fixed twisted modules), and the sum rule holds ($256$ checks).

\emph{Limitation:} due to factorization, the script does not compute $U=\Vev^G$ and 
  test Lemma \ref{lemorbrel} for the diagonal group. The outcome is implied by the 
corresponding statements for the factors ($A$ with $\Z/4$, all modules stable, and $B$ with 
$\Z/2$). The example tests the formulation for composite order with proper stabilizers, i.e., 
 index set $\Q_G$, orbit conventions, $J$-fixed twisted modules 
in the composite-order setting.  
A non-product example is listed in Section \ref{secopen}.

\end{example}

\section{Affine even part and open questions} 
\label{secopen}

\subsection{Affine even part}

 Let $\Vev=L_{\mathfrak g}(k,0)$ be an affine vertex operator algebra at positive integer 
level and $J=L_{\mathfrak g}(k\omega_J)$ a simple current of order two with 
$\theta_J=-\dim_qJ$, so that $V=\Vev\oplus J$ is a simple VOSA (e.g., $\mathfrak{sl}_2$ with 
$k\equiv2\pmod4$; $\mathfrak{sl}_{2m}$, $J=L(k\omega_m)$ with $km\equiv2\pmod4$; Examples 
\ref{exB}, \ref{exE}). 
 If $\sigma$ is a diagram automorphism of order two preserving $J$,  
then by \cite{DL} the orbifold $S$-matrix is given by the Kac-Peterson coefficients and 
the fusion coefficients $N^{\lambda_3^\dagger}_{\lambda_1,\lambda_2^\dagger}$ ($\lambda_1$ 
untwisted, $\lambda_2^\dagger$, $\lambda_3^\dagger$ twisted) 
 by the twisted Kac-Walton formula. 
 Theorem \ref{thmmain} (a) then expresses 
$N^\pm=N^{\lambda_3^\dagger}_{\lambda_1,\lambda_2^\dagger}\pm 
N^{J\fus\lambda_3^\dagger}_{\lambda_1,\lambda_2^\dagger}$ through these coefficients, 
with the parity selection rule. 
   In the affine cases of Examples \ref{exB}, \ref{exC}, \ref{exE} everything is computed 
from the modular data of $U$.  
  For inner automorphisms of affine $\Vev$ all twisted modules are stable under 
  the property of Section \ref{secinner} (Proposition \ref{propinner} (ii)).  
  Thus Theorem \ref{thmTV} (ii) and Theorem 
\ref{thmmain} (a) apply to every generating triple. The resulting coefficients agree with  
the untwisted ones when the first argument is untwisted (Proposition \ref{propinner} (i)). 

\subsection{Open problems} 

In this Subsection we list the remaining open problems related to the material of this paper. 

\emph{1. Composite order with proper stabilizers, beyond products:} example \ref{exF} 
is a product and tests only statements. A non-product test in which 
$U=\Vev^G$ is computed is open: e.g., $V=F^{\otimes4}$ with $\sigma$ a $4$-cycle 
(thus $\Vev=L_1(\mathfrak{sl}_2)^{\otimes2}$, 
 $\sigma=\mathrm{swap}\circ(\text{inner involution})$, $G=\Z/4$, the orbit 
 $\{(0,\frac12),(\frac12,0)\}$ with stabilizer $\Z/2$). 
Another open question is represente by two unstable modules in a composite-order triple 
beyond the orbit average (see Remark \ref{remTVopen}). 

\emph{2. Non-cyclic abelian and non-abelian $G$ (projective representations of stabilizers).} 

\emph{3. The correspondence for non-inner automorphisms (Example \ref{exB}):}
 existence of a tensor product of twisted $V$-modules and monoidality of $F_V$. 
 For inner automorphisms it depends on conditions \ref{hypL} (Li results for superalgebras). 

\emph{4. Type-$Q$ structure of $F_V(X)$ for $J$-fixed $X$:}
Theorem \ref{thmmain} (a) counts correctly but nothing is known beyond 
conditions \ref{hypdict}.   

\emph{5. Parity type of fixed twisted modules}: for $n$ odd, is a $J$-fixed twisted module 
 necessarily parity-twisted? See Remark \ref{remweights}.  

\emph{6. Examples with $d=-1$ and non-inner automorphisms:} Example \ref{exC} relies on 
 quoted input (a). The Kac-Walton route for $k>1$ is untested. 

\emph{7. Twisted$\times$twisted transport of the third condition of \ref{hypL}
 and the $G$-crossed identification:} 
 whether Theorem \ref{thmTV} is the cyclic case of the $G$-crossed Verlinde formula 
 \cite{T, DeM}, and a comparison of Lemma \ref{lemorbrel} with the relations of 
 \cite{DRX1} for general stabilizers.

The material of this paper is also useful in other areas 
of mathematical physics \cite{DKFRT, FZ, FFRS},   
in particular, in computations of topological invariants of condensed matter theory 
\cite{SK, BVZ, SuZ, ZA, kmmzz},  
and foliation theory \cite{Z1, Z2, Z3}. 

\section*{Acknowledgments}

The author is supported by the Institute of Mathematics, Academy of Sciences of the Czech
Republic (RVO 67985840). The author would like to thank C. Dong and S. Kanade 
for previous discussions. 

	\medskip
	\noindent\textbf{Data Availability.}
	Data sharing is not applicable to this article as no datasets were generated
	or analysed during the current study.

	\medskip
	\noindent\textbf{Declarations}

This paper arose in partial collaboration with ChatGPT, an AI assistant developed by 
OpenAI. The collaboration involved exploring proof strategies, refining arguments, locating
references, Pyton computations, and drafting. Revision and writing of final text have been done by the authors.
The authors take full responsibility for the mathematical content and the final text.
	
\medskip
	\noindent\textbf{Conflict of interest.}
	The authors have no conflicts of interest to declare that are relevant to the
	content of this article.


\end{document}